\documentclass{article}

\title{Approximating the distribution of the $L_q$-norm of a random point in a $d$-dimensional cube}
\author{Zoe Shapcott}
\date{School of Mathematics, Cardiff University, email: {\scriptsize ShapcottZA@cardiff.ac.uk} }

\usepackage[margin=2cm, top=2cm, bottom=2cm]{geometry}
\usepackage[tbtags]{amsmath}
\usepackage{amssymb}
\usepackage{amsthm}
\usepackage{graphicx}
\usepackage{listings}
\usepackage{xcolor}
\usepackage[sorting=none]{biblatex}
\usepackage{caption}
\usepackage{subcaption}
\def\dd{\mathrm{d}}

\numberwithin{figure}{section}
\numberwithin{equation}{section}

\usepackage[many]{tcolorbox}
\definecolor{main}{HTML}{000000}
\definecolor{sub}{HTML}{ffffff}
\definecolor{border}{HTML}{0000ff}
\definecolor{codegreen}{rgb}{0,0.6,0}

\definecolor{boxBborder}{HTML}{d8d8d8}
\newtcolorbox{BoxA}{
    fontupper = \color{main},
    boxrule = 1.5pt,
    colframe = border,
    colback = sub, 
    rounded corners,
    arc = 5pt
}
\newtcolorbox{BoxB}{
    fontupper = \color{main},
    boxrule = 1.5pt,
    colframe = boxBborder,
    colback = sub, 
    rounded corners,
    arc = 5pt
}

\hypersetup{colorlinks = true,linkcolor = blue,urlcolor  = blue,citecolor = blue,anchorcolor = blue,breaklinks = true}

\begin{document}

\maketitle

\begin{abstract}
    In this note, we assess the accuracy of CLT-based approximations for the volume of intersection of the \(d\)-dimensional cube \([-1,1]^{d}\) and an \(L_{q}\)-ball centred at the origin; this is clearly equivalent to approximating the distribution of the \(L_{q}\)-norm of a random point in a \(d\)-dimensional cube centered at 0. The approximations are CLT-based where to improve the normal approximation we use the first term in the Edgeworth expansion.
     We have included a section analysing the information obtained from ChatGPT in response to prompts regarding this theory; in our case, ChatGPT answers were not very helpful. Illustrations of the approximation formulae, as the radius of the ball increases, for different values of \(d\) and \(q\) are also given, alongside lines showing a Monte Carlo simulation of the intersection volume. 
\end{abstract}

\section{Introduction}\label{sec:Introduction}
In this note, we explore two formulae for approximating the volume of intersection of a \(d\)-dimensional cube and an \(L_{q}\)-ball using the Normal distribution. 
We start in Section \ref{sec:Theory} by detailing the necessary theory and developing the formulae. This includes highlighting some basic notations and key definitions, the Law of Large Numbers, the Central Limit Theorem (CLT), the first three cumulants and the first-order Edgeworth expansion.

Having constructed our two approximations, the next section (Section \ref{sec:ChatGPT}) provides a discussion of what ChatGPT provides in response to prompts regarding this topic. This is completed by comparing the information obtained from different questions with the theory described in Section \ref{sec:Theory} and other existing literature.

Following this, in Section \ref{sec:Simulations}, we present numerical experiments comparing the approximations developed in Section \ref{sec:Theory} with simulations of the intersection volume obtained via Monte Carlo methods. This is performed for different values of \(q\) for the \(L_{q}\)-ball and different dimensions \(d\). The graphs are plotted over increasing values for the radius. We are thus able to examine the effects of changing \(d\) and \(q\) on the accuracy of the approximations, with particular attention also paid to the tails of the distributions. Note that the programming language R is used to run the simulations and produce the graphs. Examples of the code used can be seen in Appendix \ref{sec:Appendix}.

\subsection{Notation}\label{subsec:Notation}
Throughout the report, let \(d > 0\) and integer denote the dimension of the space we are working in. Additionally, let a \(d\)-dimensional vector \(\mathbf{x}\) be given by
\begin{equation}
    \mathbf{x} = (x_{1},\dotsc,x_{d}).
\end{equation}
Then, for \(q > 0\), the \(L_{q}\)-norm of \(\mathbf{x}\) is
\begin{equation}
    \|\mathbf{x}\|_{q} = \left(\sum_{i=1}^{d}|x_{i}|^{q}\right)^{1/q}.
    \label{eq:LqNorm}
\end{equation}
Taking the \(q\)-th power of both sides of equation (\ref{eq:LqNorm}), we obtain
\begin{equation}
    \|\mathbf{x}\|_{q}^{q} = \sum_{i=1}^{d}|x_{i}|^{q}.
    \label{eq:qthPowerLqNorm}
\end{equation}
This represents the \(q\)-th power of the distance between \(\mathbf{x}\) and the origin, \(\mathbf{0}=(0,\dotsc,0)\), and will be an important quantity in the analysis of Section \ref{sec:Theory}.

\section{Theory}\label{sec:Theory}
In this section we aim to develop two approximations for the volume of intersection of a cube and an \(L_{q}\)-ball. These are a Normal approximation based on the Central Limit Theorem (Section \ref{subsec:NormalApprox}) and an improved Normal approximation using the first order Edgeworth Expansion (Section \ref{subsec:ImprovedNormalApprox}). Specifically, we are interested in approximating the volume of intersection between the \(d\)-dimensional cube \([-1,1]^{d}\) and the \(L_{q}\)-ball of radius \(r>0\) centred at the origin. We denote this ball by \(\mathcal{B}_{d,q}(r)\) and define it via
\begin{equation}
    \mathcal{B}_{d,q}(r) = \left\{\mathbf{x} \in \mathbb{R}^d : \|\mathbf{x}\|_{q} \leq r \right\},
    \label{eq:LqBall}
\end{equation}
where \(\|\mathbf{x}\|_{q}\) is as given in equation (\ref{eq:LqNorm}).

In order to estimate the volume of intersection of the cube and the \(L_{q}\)-ball, we look for a formula to approximate the fraction of the cube covered by the ball. This is completed in a similar way to the specific case of \(q=2\) considered in detail in \cite{zhigljavsky2020covering} and 
\cite{noonan2023improving}; see also  Section \(1.1.8\) of \cite{AnatolyBook}. We define this quantity as
\begin{equation}
    C_{d,q,r} = \frac{\text{vol}([-1,1]^{d} \cap \mathcal{B}_{d,q}(r))}{2^{d}},
\end{equation}
where \(\text{vol}(A)\) denotes the volume of \(A\). The denominator is \(2^{d}\) since \(\text{vol}([-1,1]^{d}) = \left(1-(-1)\right)^{d} = 2^{d}\). Note that \(C_{d,q,r}\) depends on the dimension \(d\), the \(L_{q}\)-norm used to define the ball, and the radius \(r\) of the ball. We may also observe that \(0 \leq C_{d,q,r} \leq 1\).

Let \(\mathbf{x}\) be a random vector uniformly distributed on the cube \([-1,1]^d\). Then we can also define \(C_{d,q,r}\) as
\begin{equation}
    C_{d,q,r} = \mathbb{P}\left\{\|\mathbf{x}\|_{q} \leq r \right\} = \mathbb{P}\left\{\|\mathbf{x}\|_{q}^{q} \leq r^{q} \right\} = \mathbb{P}\left\{\sum_{i=1}^{d}|x_{i}|^{q} \leq r^{q} \right\}.
    \label{eq:CdqrX}
\end{equation}
Thus to estimate the fraction of the cube \([-1,1]^{d}\) covered by the \(L_{q}\)-ball, we can try to approximate the cumulative distribution function (CDF) of \(\sum_{i=1}^{d}|x_{i}|^{q}\).

Since \(\mathbf{x}=(x_{1},\dotsc,x_{d})\) is uniform random on \([-1,1]^{d}\), we have that the \(x_{i}\) (\(i = 1,\dotsc,d\)) are all independent and identically distributed random variables. Specifically, each \(x_{i}\) is uniformly distributed on \([-1,1]\), where the probability density function (PDF) of each \(x_{i}\) is given by \(f_{x}(t) = \frac{1}{2}\) for \(t \in [-1,1]\) and \(f_{x}(t)=0\) otherwise. Consequently, we also have that \(\sum_{i=1}^{d}|x_{i}^{q}|\) is equal in distribution to \(\sum_{i=1}^{d}z_{i}^{q}\), where the \(z_{i}\) \((i = 1,\dotsc,d)\) are independent and identically distributed random variables, uniform on \([0,1]\). This gives us that the PDF of each \(z_{i}\) is \(f_{z}(t)=1\) for \(t \in [0,1]\) and zero otherwise. Using this and equation (\ref{eq:CdqrX}), then yields
\begin{equation}
    C_{d,q,r} = \mathbb{P}\left\{\sum_{i=1}^{d}z_{i}^{q} \leq r^{q} \right\},
    \label{eq:Cdqrz}
\end{equation}
and to derive the desired formulae we search for methods of approximating the CDF of \(\sum_{i=1}^{d}z_{i}^{q}\). We begin by looking for an approximation using the Normal distribution.

\subsection{Normal Approximation for \texorpdfstring{\(C_{d,q,r}\)}{Cd,q,r}}\label{subsec:NormalApprox}
We break down finding a Normal approximation to the distribution of \(\sum_{i=1}^{d}z_{i}^{q}\), where the \(z_{i}\) are independent random variables with uniform distribution on \([0,1]\) (as defined above in Section \ref{sec:Theory}), into a few steps. First, in Section \ref{subsubsec:Moments}, we calculate the moments of order \(\alpha\) of the random variables \(z_{i}^{q}\). Subsequently, in Section \ref{subsubsec:LLN}, we apply the Law of Large Numbers to \(\frac{1}{d}\sum_{i=1}^{d}z_{i}^{q}\). Finally, in Section \ref{subsubsec:CLT}, we derive the desired formula using the Central Limit Theorem.

\subsubsection{Moments of \texorpdfstring{\(z_{i}^{q}\)}{z\_i\^\ q}}\label{subsubsec:Moments}
In order to find a Normal approximation for the distribution of \(\sum_{i=1}^{d}z_{i}^{q}\) we first need to calculate the raw moments of order \(\alpha\) of \(z_{i}^{q}\).

Recall that for a continuous random variable \(X\) with PDF \(f_{X}(x)\) on \([a,b]\), the raw moments, or moments about the origin, of order \(\alpha\) are defined as
\begin{equation}
    \mu_{\alpha} = \mathbb{E}(X^{\alpha}) = \int_{a}^{b}x^{\alpha}f_{X}(x)\dd x,
    \label{eq:MomentsDef}
\end{equation}
see Section 2.5 of \cite{IntroductionToProbabilityBook} for more detail. Therefore, since the PDF of each \(z_{i}\) \((i = 1,\dotsc,d)\) is given by \(f_{z}(t) = 1\) for \(t \in [0,1]\) and zero otherwise (see Section \ref{sec:Theory}), we can calculate that the \(z_{i}^{q}\) raw moments of order \(\alpha\) are given by
\begin{equation}
    \mu_{\alpha,q} = \frac{1}{\alpha q + 1},
    \label{eq:Moments}
\end{equation}
for \(\alpha > 0\) and \(q > 0\). Note that equation (\ref{eq:Moments}) follows from direct computation using the definition in equation (\ref{eq:MomentsDef}), as
\begin{equation}
    \begin{split}
        \mu_{\alpha,q} & = \mathbb{E}\left[(z_{i}^{q})^{\alpha}\right] = \mathbb{E}\left[z_{i}^{\alpha q}\right] \\
        & = \int_{0}^{1}x^{\alpha q} \cdot 1 \dd x \\
        & = \left[\frac{x^{\alpha q +1}}{\alpha q + 1}\right]_{0}^{1} \\
        & = \frac{1}{\alpha q + 1}. \\
    \end{split}
\end{equation}

In particular, we can use equation (\ref{eq:Moments}) to find that the mean of \(z_{i}^{q}\) is given by
\begin{equation}
    \mathbb{E}(z_{i}^{q}) = \mu_{1,q} = \frac{1}{q+1}.
    \label{eq:ziqMean}
\end{equation}
Additionally, we have
\begin{equation}
    \mathbb{E}\left[(z_{i}^{q})^{2}\right] = \mu_{2,q} = \frac{1}{2q+1}.
\end{equation}
The variance of \(z_{i}^{q}\) can then be calculated via
\begin{equation}
    \begin{split}
    \mathbb{V}\text{ar}(z_{i}^{q}) & = \mathbb{E}\left[(z_{i}^{q})^{2}\right] - \left[\mathbb{E}(z_{i}^{q})\right]^{2} \\
    & = \mu_{2,q} - \mu_{1,q}^{2} \\
    & = \frac{1}{2q+1} - \frac{1}{(q+1)^{2}}. \\
    \end{split}
    \label{eq:VarStep}
\end{equation}
Simplifying the above, we thus obtain
\begin{equation}
    \mathbb{V}\text{ar}(z_{i}^{q}) = \frac{q^{2}}{(2q+1)(q+1)^{2}}.
    \label{eq:ziqVar}
\end{equation}

Using equation (\ref{eq:ziqMean}) for the expectation and equation (\ref{eq:ziqVar}) for the variance of the random variables \(z_{i}^{q}\), we next apply the Law of Large Numbers.

\subsubsection{The Law of Large Numbers}\label{subsubsec:LLN}
In this section, We first recall the Strong Law of Large Numbers using Section \(8.3\) of \cite{IntroductionToProbabilityBook}, and then apply it to our specific case.

The Law of Large Numbers states that, for a sequence of independent and identically distributed random variables \(X_{1}, X_{2}, \dotsc\, X_{n}\), with finite mean \(\mu\) and finite variance \(\sigma^{2}\), we have that
\begin{equation}
    \overline{X}_{n} = \frac{1}{n}(X_{1} + X_{2} + \cdots + X_{n}) \longrightarrow \mu,
\end{equation}
with probability one as \(n \longrightarrow \infty\).

Implementing this in our case, we see that
\begin{equation}
    \frac{1}{d}(z_{1}^{q} + z_{2}^{q} + \cdots + z_{d}^{q}) \longrightarrow \mathbb{E}(z_{i}^{q}),
\end{equation}
with probability one as \(d \longrightarrow \infty\). This is because the \(z_{i}^{q}\) are all independent and identically distributed random variables with finite mean and finite variance, as computed in equations (\ref{eq:ziqMean}) and (\ref{eq:ziqVar}). Hence by the Law of Large Numbers, we obtain that
\begin{equation}
    \frac{1}{d}\sum_{i=1}^{d}z_{i}^{q} \longrightarrow \frac{1}{q+1},
\end{equation}
with probability one as \(d \longrightarrow \infty\).

In the next section we use this with the Central Limit Theorem to derive the desired approximation formula.

\subsubsection{The Central Limit Theorem}\label{subsubsec:CLT}
Within this section, we begin by recalling the Central Limit Theorem, before utilising this to derive our first approximation formula for \(C_{d,q,r}\), the fraction of the cube covered by the \(L_{q}\)-ball, as defined in equation (\ref{eq:Cdqrz}).

Let \(X_{1},X_{2},\dotsc,X_{n}\) be a sequence of independent and identically distributed random variables with finite mean \(\mu\) and finite nonzero variance \(\sigma^{2}\). By the Law of Large Numbers (Section \ref{subsubsec:LLN}) we have that \(\overline{X}_{n} = \frac{1}{n}\sum_{i=1}^{n}X_{i} \longrightarrow \mu\) as \(n \longrightarrow \infty\) with probability one. Now, consider
\begin{equation}
    \frac{S_{n} - \mathbb{E}(S_{n})}{\sqrt{\mathbb{V}\text{ar}(S_{n})}},
\end{equation}
where \(S_{n} = n\overline{X}_{n} = \sum_{i=1}^{n}X_{i}\), see Section 8.3 of \cite{ProbabilityAnIntroductionBook}. We have
\begin{equation}
     \mathbb{E}(S_{n}) =  \mathbb{E}\left(\sum_{i=1}^{n}X_{i}\right) = \sum_{i=1}^{n}\mathbb{E}(X_{i}) = n\mu,
\end{equation}
and, since the \(X_{i}\) are all independent,
\begin{equation}
    \mathbb{V}\text{ar}(S_{n}) = \mathbb{V}\text{ar}\left(\sum_{i=1}^{n}X_{i}\right) = \sum_{i=1}^{n}\mathbb{V}\text{ar}(X_{i}) = n\sigma^{2}.
\end{equation}
Thus
\begin{equation}
    \frac{S_{n} - \mathbb{E}(S_{n})}{\sqrt{\mathbb{V}\text{ar}(S_{n})}} = \frac{S_{n}-n\mu}{\sigma\sqrt{n}}.
\end{equation}
The Central Limit Theorem then states that
\begin{equation}
    \mathbb{P}\left\{\frac{S_{n}-n\mu}{\sigma\sqrt{n}} \leq t\right\} \longrightarrow \Phi(t),
    \label{eq:CLTSn}
\end{equation}
as \(n \longrightarrow \infty\) for all \(t \in \mathbb{R}\), where
\begin{equation}
    \Phi(t) = \frac{1}{\sqrt{2\pi}}\int_{-\infty}^{t}e^{-\frac{1}{2}u^{2}} \dd u,
\end{equation}
is the CDF of the Standard Normal distribution, \(\mathcal{N}(0,1)\). See Section 8.3 of \cite{ProbabilityAnIntroductionBook} for more details.

Rearranging the inequality in equation (\ref{eq:CLTSn}), we can see that it is equivalent to say that
\begin{equation}
    \mathbb{P}\left\{S_{n} \leq t\sigma\sqrt{n} + n\mu\right\} \longrightarrow \Phi(t),
\end{equation}
as \(n \longrightarrow\infty\). Therefore, using some simple algebraic manipulations, we obtain that for a constant \(c\),
\begin{equation}
    \mathbb{P}\left\{S_{n} \leq c\right\} \longrightarrow \Phi\left(\frac{c - n\mu}{\sigma\sqrt{n}}\right),
    \label{eq:FirstApproxGeneral}
\end{equation}
as \(n \longrightarrow \infty\).

We are now in a position to obtain our first approximation by substituting the relevant values into equation (\ref{eq:CLTSn}) or (\ref{eq:FirstApproxGeneral}). Recall that we have a sequence of \(d\) independent and identically distributed random variables \(z_{1}^{q}, z_{2}^{q}, \dotsc, z_{d}^{q} \), where each \(z_{i}\) is uniformly distributed on \([0,1]\). These random variables \(z_{i}^{q}\) have finite mean \(\mu\) and finite nonzero variance \(\sigma^{2}\) given respectively by
\begin{equation}
    \mu = \frac{1}{q+1} \quad \text{ and } \quad \sigma^{2} = \frac{q^{2}}{(2q+1)(q+1)^{2}},
\end{equation}
as found in equations (\ref{eq:ziqMean}) and (\ref{eq:ziqVar}). Therefore, by the Central Limit Theorem as stated in equation (\ref{eq:CLTSn}), we obtain
\begin{equation}
    \mathbb{P}\left\{\frac{\sum_{i=1}^{d}z_{i}^{q} - \frac{d}{q+1}}{\sqrt{d}\sqrt{\frac{q^{2}}{(2q+1)(q+1)^{2}}}} \leq t\right\} \longrightarrow \Phi(t),
\end{equation}
as \(d \longrightarrow \infty\). This can be alternatively stated using equation (\ref{eq:FirstApproxGeneral}) as
\begin{equation}
    \mathbb{P}\left\{\sum_{i=1}^{d}z_{i}^{q} \leq t\right\} \longrightarrow \Phi\left(\frac{t-\frac{d}{q+1}}{\sqrt{d}\sqrt{\frac{q^{2}}{(2q+1)(q+1)^{2}}}}\right),
    \label{eq:Approx1}
\end{equation}
for \(d \longrightarrow \infty\).

We can simplify the argument of \(\Phi\) in equation (\ref{eq:Approx1}) and let \(t = r^{q}\) to obtain our first approximation formula for \(C_{d,q,r}\). We investigate the effectiveness of this approximation for different values of \(d\) and \(q\) in Section \ref{sec:Simulations}.
\begin{BoxA}
    \textbf{Approximation 1: Normal Approximation}

    For \(d\) large enough, we can approximate \(C_{d,q,r}\), by
    \begin{equation}
         C_{d,q,r} = \mathbb{P}\left\{\|\mathbf{x}\|_{q}^{q} \leq r^{q}\right\} = \mathbb{P}\left\{\sum_{i=1}^{d}z_{i}^{q} \leq r^{q}\right\} \approx \Phi\left(\frac{\sqrt{2q+1}\bigl((q+1)r^{q}- d\bigr)}{q\sqrt{d}}\right).
        \label{eq:Approx1BoxVersion}
    \end{equation}
\end{BoxA}

Having now established the first of our approximations for \(C_{d,q,r}\) using the normal distribution, in Section \ref{subsec:ImprovedNormalApprox} we next look to find an improvement on this formula using the first order Edgeworth Expansion.

\subsection{Improved Normal Approximation for \texorpdfstring{\(C_{d,q,r}\)}{Cd,q,r}}\label{subsec:ImprovedNormalApprox}
In this section, similarly to Section \ref{subsec:NormalApprox}, we split the process of deriving a second, more accurate approximation for the distribution of \(\sum_{i=1}^{d}z_{i}^{q}\) into a couple of stages. Note that the \(z_{i}\) are independent and identically distributed random variables with uniform distribution on \([0,1]\), as described in Section \ref{sec:Theory}. Firstly, in Section \ref{subsubsec:Cumulants}, we define and calculate the first few cumulants of \(z_{i}^{q}\) (\(i=1,\dotsc,d\)). Following this, in Section \ref{subsubsec:Edgeworth}, we introduce the first order Edgeworth Expansion, which is an extension of the Central Limit Theorem. We then deduce a second formula approximating the fraction of the cube covered by the \(L_{q}\)-ball by adding this term in the expansion to the formula (\ref{eq:Approx1BoxVersion}) found in Section \ref{subsubsec:CLT}. The relative accuracy of the two formulas against a Monte Carlo simulation of the quantity will be explored in Section \ref{sec:Simulations}.

\subsubsection{Cumulants of \texorpdfstring{\(z_{i}^{q}\)}{z\_i\^\ q}}\label{subsubsec:Cumulants}
In order to apply the first order Edgeworth Expansion in Section \ref{subsubsec:Edgeworth}, we need the first three cumulants of the random variables \(z_{i}^{q}\) defined in Section \ref{sec:Theory}. In general, cumulants of order \(n\) can be found using a cumulant generating function (see \cite{CumulantWolfram} for more details) but the form of the cumulant can be quite complicated for general n. However, here we are only interested in the first three, which can be simply defined as follows.

For a random variable \(\xi\), the first three cumulants, \(\text{cum}_{k}(\xi)\), \(k = 1,2,3\), are given by
\begin{align}
    \text{cum}_{1}(\xi) & = \mathbb{E}(\xi) \qquad \qquad \qquad (\text{the mean of } \xi), \\
    \text{cum}_{2}(\xi) & = \mathbb{V}\text{ar}(\xi) \qquad \qquad \quad (\text{the variance of } \xi), \text{ and }\\
    \text{cum}_{3}(\xi) & = \mathbb{E}\left[\Bigl(\xi - \mathbb{E}(\xi)\Bigr)^{3}\right] \; \; (\text{the third central moment of } \xi).
\end{align}

From the results of equations (\ref{eq:ziqMean}) and (\ref{eq:ziqVar}) in Section \ref{subsubsec:Moments}, we already know that, for all \(i=1,\dotsc,d\), the first cumulant of \(z_{i}^{q}\) is
\begin{equation}
    \text{cum}_{1}(z_{i}^{q}) = \mathbb{E}(z_{i}^{q}) = \frac{1}{q+1},
\end{equation}
and the second cumulant of \(z_{i}^{q}\) is
\begin{equation}
    \text{cum}_{2}(z_{i}^{q}) = \mathbb{V}\text{ar}(z_{i}^{q}) = \frac{q^{2}}{(2q+1)(q+1)^{2}}.
\end{equation}
Thus it remains to calculate the third cumulant of \(z_{i}^{q}\). Using the definition, We have
\begin{equation}
    \begin{split}
        \text{cum}_{3}(z_{i}^{q}) & = \mathbb{E}\left[\bigl(z_{i}^{q} - \mathbb{E}(z_{i}^{q})\bigr)^{3}\right] = \mathbb{E}\left[\left(z_{i}^{q} - \frac{1}{q+1}\right)^{3}\right] \\
        & = \int_{0}^{1}\left(t^{q} - \frac{1}{q+1}\right)^{3} \cdot 1 \dd t, \\
        \label{eq:CumulantSubstep}
    \end{split}
\end{equation}
since the PDF of \(z_{i}\) is equal to one for all t in the interval \([0,1]\) and zero otherwise. Expanding the right hand side of equation (\ref{eq:CumulantSubstep}), we get
\begin{equation}
    \begin{split}
        \text{cum}_{3}(z_{i}^{q}) & = \int_{0}^{1} \left[t^{3q} - \frac{3t^{2q}}{q+1} + \frac{3t^{q}}{(q+1)^{2}} - \frac{1}{(q+1)^{3}}\right] \dd t \\
        & = \frac{1}{3q+1} - \frac{3}{(2q+1)(q+1)} + \frac{2}{(q+1)^{3}}. \\
    \end{split}
\end{equation}
Performing some basic algebraic manipulations on this expression yields
\begin{equation}
    \text{cum}_{3}(z_{i}^{q}) = \frac{2q^{3}(q-1)}{(3q+1)(2q+1)(q+1)^{3}},
    \label{eq:ThirdCumulant}
\end{equation}
and we are now ready to introduce the first order Edgeworth Expansion to approximate the distribution of \(\|\mathbf{x}\|_{q}^{q}\).

\subsubsection{The Edgeworth Expansion}\label{subsubsec:Edgeworth}
To develop an improved version of our first formula estimating \(C_{d,q,r}\), equation (\ref{eq:Approx1BoxVersion}), we include the first term in the Edgeworth Expansion, which is made up of correction terms added to the Central Limit Theorem. For a general discussion of all terms in the expansion, see  for example Theorem \(7.9\) in \cite{IntroductionToStatisticalLimitTheory}. Here, however, we use Theorem 7.8 in \cite{IntroductionToStatisticalLimitTheory} to define the first order Edgeworth Expansion. For \(X_{1},X_{2},\dotsc,X_{n}\) a sequence of \(n\) independent and identically distributed random variables with mean \(\mu\) and variance \(\sigma^{2}\), the first order Edgeworth Expansion states that
\begin{equation}
    \mathbb{P}\left\{\frac{\sum_{i=1}^{n}X_{i} - n\mu}{\sigma\sqrt{n}} \leq t\right\} \longrightarrow \Phi(t) + \frac{(1-t^{2})\phi(t)}{6 \sqrt{n}} \cdot \frac{\text{cum}_{3}(X_{i})}{\sigma^{3}},
\end{equation}
as \(n \longrightarrow \infty\). Here \(\phi(t)\) is the PDF of the standard Normal distribution, \(\mathcal{N}(0,1)\), and it is defined via
\begin{equation}
    \phi(t) = \frac{1}{2\pi}e^{-t^{2}/2}.
\end{equation}

We now apply this formula to our situation. The \(z_{i}^{q}, \; (i=1,2,\dotsc,d),\) as defined in Section \ref{sec:Theory} are all independent and identically distributed random variables. These have mean \(\mu = 1/(q+1)\) from equation (\ref{eq:ziqMean}), variance \(\sigma^{2} = q^{2}/(2q+1)(q+1)^{2}\) from equation (\ref{eq:ziqVar}), and third cumulant \(\text{cum}_{3}(z_{i}^{q}) = 2q^{3}(q-1)/(3q+1)(2q+1)(q+1)^{3}\) from equation (\ref{eq:ThirdCumulant}). We thus obtain that the first order Edgeworth Expansion implies
\begin{equation}
    \mathbb{P}\left\{\frac{\sum_{i=1}^{d}z_{i}^{q} - \frac{d}{q+1}}{\sqrt{d}\sqrt{\frac{q^{2}}{(2q+1)(q+1)^{2}}}} \leq t\right\} \longrightarrow \Phi(t) + \frac{(1-t^{2})\phi(t)}{6\sqrt{d}} \cdot \frac{2q^{2}(q-1)}{(3q+1)(2q+1)(q+1)^{3}} \cdot \left(\sqrt{\frac{(2q+1)(q+1)^{2}}{q^{2}}}\right)^{3},
    \label{eq:EdgeworthStep}
\end{equation}
as \(d \longrightarrow \infty\).

Rearranging and simplifying the inequality on the left of equation (\ref{eq:EdgeworthStep}) and also simplifying the expression on the right yields that
\begin{equation}
    \mathbb{P}\left\{\sum_{i=1}^{d}z_{i}^{q} \leq \frac{t\sqrt{d}q}{(q+1)\sqrt{2q+1}} + \frac{d}{q+1}\right\} \longrightarrow \Phi(t) + \frac{(1-t^{2})\phi(t)}{6\sqrt{d}} \cdot \frac{2(q-1)\sqrt{2q+1}}{3q+1},
\end{equation}
as \(d \longrightarrow \infty\).

Therefore, we can now state our second Normal distribution based approximation for \(C_{d,q,r}\). We compare the results of using this formula to the one obtained in equation (\ref{eq:Approx1BoxVersion}) in the graphs and discussions of Section \ref{sec:Simulations}.
\begin{BoxA}
    \textbf{Approximation 2: Improved Normal Approximation}

    For \(d\) large enough, we can approximate \(C_{d,q,r}\), by
    \begin{equation}
         C_{d,q,r} = \mathbb{P}\left\{\|\mathbf{x}\|_{q}^{q} \leq r^{q}\right\} = \mathbb{P}\left\{\sum_{i=1}^{d}z_{i}^{q} \leq r^{q}\right\} \approx \Phi(t_{r}) + \frac{(1-t_{r}^2)\phi(t_{r})}{6\sqrt{d}} \cdot \frac{2(q-1)\sqrt{2q+1}}{3q+1},
        \label{eq:Approx2BoxVersion}
    \end{equation}
    where
    \begin{equation}
        t_{r} = \frac{\sqrt{2q+1}\bigl((q+1)r^{q} - d\bigr)}{q\sqrt{d}}.
        \label{eq:Approx2BoxVersiontrValue}
    \end{equation}
\end{BoxA}

Now that we have developed two formulae to approximate \(C_{d,q,r}\) (the fraction of the cube covered by the \(L_{q}\)-ball) as a function of the radius \(r\) for different values of \(d\) and \(q\), we next compare the results found here with information obtained from ChatGPT (Section \ref{sec:ChatGPT}). Following this, in Section \ref{sec:Simulations}
we will graph both these formulae against Monte Carlo simulations of the quantity \(C_{d,q,r}\) as \(r\) increases for a variety of different values for \(d\) and \(q\).

\section{ChatGPT}\label{sec:ChatGPT}
In this section we analyse information obtained from ChatGPT in response to prompts linked to the theory discussed in Section \ref{sec:Theory}, in order to assess the accuracy and scope of knowledge it has with regards to this topic.

The first question asked to ChatGPT investigates what it could say about the generic topic of exploration in this report. The specific prompt used was
\begin{BoxB}
    {\small``What is the formula for the volume of intersection of the d-dimensional cube [-1,1]\textasciicircum{}d and an L\_q-ball of radius r centred at the origin.''}
\end{BoxB}
which resulted in the response given below \cite{ChatGPT}.
\begin{BoxB}
    {\small``The problem you're describing involves finding the volume of the intersection between a d-dimensional cube and an L\_q-ball of radius r centered at the origin. This is a complex mathematical problem, and the formula for the volume of such an intersection may not have a simple closed form, especially as you increase the dimensionality.

    The notation "L\_q-ball" refers to a norm ball in L\_q space, where the norm is defined by the L\_q norm. The L\_q norm of a vector x = (x\_1, x\_2, ..., x\_d) in d-dimensional space is given by:

    \[ \lVert x \rVert_q = \left( \lvert x_1 \rvert^q + \lvert x_2 \rvert^q + \ldots + \lvert x_d \rvert^q \right)^{1/q} \]

    The intersection you're interested in is the common region between the d-dimensional cube \([-1, 1]^d\) and the L\_q-ball of radius r.

    Unfortunately, finding a closed-form expression for the volume of this intersection is challenging, especially as the dimensionality increases. In such cases, numerical methods or specialized software may be used to approximate the volume.

    If you have a specific value for \( q \) and \( r \), you may be able to find more information or computational approaches tailored to your exact problem in research literature or specialized mathematical software documentation.'' }
\end{BoxB}
From the above response, it may appear that ChatGPT \cite{ChatGPT} has, in general, understood the intention of the question and it has recognised that it is a complex problem depending on the values of \(q, r\) and \(d\). However, it has not been able to provide the asked for formula, nor an approximation that might be appropriate, instead claiming that numerical simulations should be used. Conversely, it has been able to correctly define the \(L_q\)-norm and the intersection that we are interested in.

When the question was simplified to the case of a Euclidean ball (i.e. \(q=2\)), via the prompt
\begin{BoxB}
    {\small``What is the formula for the volume of intersection of the d-dimensional cube [-1,1]\textasciicircum{}d and the Euclidean ball of radius r centred at the origin?''}
\end{BoxB}
a much more detailed response was obtained \cite{ChatGPT}. This may suggest that the case of \(q=2\) is much more widely studied and, as a result, more information is available for this scenario. The response obtained is shown in the next box below \cite{ChatGPT}. By comparing the formula presented in this answer generated by ChatGPT with the one documented in Section \(1.1.1\) of \cite{AnatolyBook}, we can see that the system has generated an accurate formula for the volume of a \(d\)-dimensional Euclidean ball in this case. However, the formula for the intersection volume is not true in general, only holding for small values of \(r\), which ChatGPT has not acknowledged.
\begin{BoxB}
    {\small``The volume of the intersection between a d-dimensional cube \([-1,1]^d\) and a Euclidean ball (also known as a hypersphere) of radius \(r\) centered at the origin can be expressed as a formula. This is often referred to as the volume of the intersection between a cube and a ball.

    For the case of a d-dimensional cube and a Euclidean ball, the formula for the volume of the intersection is given by:

    \[ V_d(r) = 2^d \cdot V_{\text{ball}}(r) \]

    where \(V_{\text{ball}}(r)\) is the volume of the d-dimensional Euclidean ball of radius \(r\), and \(2^d\) accounts for the symmetry due to the cube.

    The formula for the volume of a d-dimensional Euclidean ball is given by:

    \[ V_{\text{ball}}(r) = \frac{\pi^{d/2}}{\Gamma\left(\frac{d}{2} + 1\right)} \cdot r^d \]

    where \(\Gamma\) is the gamma function.

    So, combining the two formulas, you get the expression for the volume of the intersection:

    \[ V_d(r) = 2^d \cdot \frac{\pi^{d/2}}{\Gamma\left(\frac{d}{2} + 1\right)} \cdot r^d \]

    This formula provides the volume of the intersection of a d-dimensional cube and a Euclidean ball of radius \(r\) centered at the origin.''}
\end{BoxB}

Initially remaining in the case of a Euclidean ball (\(q=2\)), since ChatGPT seems to have more to say in this scenario, we now prompt ChatGPT further to try to extract the first of our approximation formulae (equation~(\ref{eq:Approx1BoxVersion})) based on the Central Limit Theorem. The question provided was
\begin{BoxB}
    {\small``What is the formula to approximate the fraction of the d-dimensional cube [-1,1]\textasciicircum{}d covered by a Euclidean ball of radius r centred at the origin using the Central Limit Theorem?''}
\end{BoxB}
and this generated the output given in the box below \cite{ChatGPT}.
\begin{BoxB}
    {\small``To approximate the fraction of the d-dimensional cube \([-1, 1]^d\) covered by a Euclidean ball of radius \(r\) centered at the origin using the Central Limit Theorem (CLT), you can use the following formula:

    \[ P(\|X\|_2 \leq r) \approx \Phi\left(\frac{r\sqrt{d}}{\sqrt{2}}\right) \]

    Here, \(X\) is a random vector representing a point chosen uniformly at random from the cube \([-1, 1]^d\), \(\|\cdot\|_2\) denotes the Euclidean norm, \(\Phi(\cdot)\) is the cumulative distribution function (CDF) of the standard normal distribution, \(r\) is the radius of the Euclidean ball, and \(d\) is the dimensionality of the space.

    This formula is based on the idea that when \(d\) is large, the distribution of the Euclidean norm of a random vector in \([-1, 1]^d\) approaches a normal distribution by the Central Limit Theorem. The factor \(r\sqrt{d}/\sqrt{2}\) comes from scaling considerations.

    Keep in mind that this is an approximation and becomes more accurate as \(d\) becomes large. Additionally, this formula assumes that the Euclidean ball is centered at the origin. If the ball is centered elsewhere, you may need to make adjustments to the formula accordingly.''}
\end{BoxB}
This response should produce results directly comparable to the formula for \(C_d,q,r\) as computed in Section \ref{subsec:NormalApprox} for the case \(q=2\) and to that of \(C_{d,Z,r}\) for \(Z=\mathbf{0}\) from Section 1.1.8 of \cite{AnatolyBook}. However, the approximation formula ChatGPT has given us is incorrect. We can see that the second sentence appears to define all the terms correctly. Additionally, ChatGPT has highlighted in the final statement the key fact that the approximation is only valid for sufficiently large \(d\) and for the specific case where the ball is centred at the origin. More importantly though, the formula it provided is wrong overall. This is because the distribution of \(\|X\|_{2}\) does not tend towards the standard Normal distribution since it does not have variance one. Thus, despite the fact that this formula would be slightly different compared to that obtained in Section \ref{subsec:NormalApprox} and \cite{AnatolyBook} due to not considering the squared norm of \(X\), the formula presented by ChatGPT is incorrect. Therefore, although it has understood the question and produced a similar result to that expected, there are serious reliability considerations of such sources, since it claims that this formula is true despite the fact that it contains errors.

Having established this, we now try to prompt ChatGPT to give us an approximation formula for general \(q\). The question
\begin{BoxB}
    {\small``What is the formula to approximate the fraction of the d-dimensional cube [-1,1]\textasciicircum{}d covered by an L\_q ball of radius r centred at the origin using the Central Limit Theorem?''}
\end{BoxB}
was asked which resulted in the output presented below \cite{ChatGPT}.
\begin{BoxB}
    {\small``The Central Limit Theorem (CLT) is typically used to approximate the distribution of sums of independent and identically distributed random variables. It's not directly used for approximating the fraction of a d-dimensional cube covered by an L\_q ball.

    However, it's possible that you're referring to a result related to the volume of a high-dimensional sphere in the context of the CLT. The volume of an \(L_q\) ball (or a sphere in this case) in d-dimensional space can be approximated using the CLT. In particular, for large \(d\), the distribution of the distance of a point chosen uniformly at random from the cube \([-1,1]^d\) to the origin converges to a normal distribution.

    One result related to this is known as the Dvoretzky's theorem, which states that for sufficiently large \(d\), most of the volume of a high-dimensional cube is concentrated in a thin shell around its boundary.

    To express this more formally, let \(B_q(r)\) be the \(L_q\) ball of radius \(r\) centered at the origin in \(d\)-dimensional space, and let \(V_d\) be the volume of the unit cube \([-1,1]^d\). Then, for large \(d\), the volume of the \(L_q\) ball \(B_q(r)\) can be approximated as:

    \[V(B_q(r)) \approx V_d \cdot P(\|X\|_q \leq r)\]

    where \(X\) is a random vector chosen uniformly at random from the cube \([-1,1]^d\) and \(\|X\|_q\) is the \(L_q\) norm of \(X\). The probability \(P(\|X\|_q \leq r)\) can be approximated using the CLT, assuming certain conditions are satisfied.

    It's important to note that the specific formula and conditions may depend on the exact context and assumptions of the problem you are considering. If you have a particular formula or context in mind, please provide more details for a more accurate and specific answer.''}
\end{BoxB}
The answer output here appears incorrect, since it claims that the Central Limit Theorem is not used in this context. However, we have utilised the Central Limit Theorem in Section \ref{subsec:NormalApprox} to derive a formula approximating the fraction of the \(d\)-dimensional cube covered by an \(L_{q}\)-ball centred at the origin. Additionally, ChatGPT instead gives information unrelated to the question asked, within which it appears to have merged together a few different concepts and not always correctly. Whilst it is true that for a high-dimensional cube, most of the volume is concentrated in an annulus near the boundary \cite{AnatolyBook}, its explanation of Dvoretzky's theorem does not completely agree with others, such as the the one in \cite{DvoretzkyWolfram}, which states that the theorem is concerned with all high-dimensional centred convex bodies and not just cubes.

Overall, it seems that ChatGPT may be able to provide some starting or basic insights in relation to this topic, particularly for the more widely studied case of \(q=2\). However, the information it relays is not fully accurate, and it may also not be clear whether or not the information is reliable. Therefore, we conclude that ChatGPT should be used with caution in this area since it is not fully informed about this specific topic of mathematics.

In the next section (Section \ref{sec:Simulations}) we present some graphical representations of the two approximation formulae and Monte Carlo simulated results for different values of \(d\) and \(q\).

\section{Numerical Simulations}\label{sec:Simulations}
In this section we present some simulations to try to assess the accuracy of the approximation formulae (\ref{eq:Approx1BoxVersion}) and (\ref{eq:Approx2BoxVersion}) for different values of \(q\) and \(d\) plotted over a range of \(r\). We consider graphs over a large interval of \(r\) to gain an appreciation of the overall distribution, with of \(C_{d,q,r}\) taking values over the whole interval \([0,1]\). Additionally, we also display some figures over a narrower range of \(r\) that specifically zoom into the tails of the distribution to acquire a deeper understanding of the reliability of the formulae there.

The code producing the graphs is written using the programming language R, and examples of the code used can be seen in Appendix \ref{sec:Appendix}. In all the graphs, the solid black line depicts the Monte Carlo simulation of the quantity of interest (i.e the fraction of the cube \([-1,1]^d\) covered by the \(L_{q}\)-ball of radius \(r\) centred at the origin, \(C_{d,q,r}\)). This gives us the closest representation of the true value for \(C_{d,q,r}\) that we can simulate, and we may expect that it increases in accuracy the more points that are simulated in the cube. The red dashed line illustrates the first Normal approximation formula as in equation (\ref{eq:Approx1BoxVersion}). Finally, the blue dot-dashed line portrays the improved Normal approximation formula, as represented in equations (\ref{eq:Approx2BoxVersion}) and (\(\ref{eq:Approx2BoxVersiontrValue}\)).

We split the investigations into a couple of sections. Firstly, in Section \ref{subsec:q2,5,10,d5,10,25,50} we compare graphs having \(q = 2, 5\) or \(10\) and \(d = 5, 10, 25\) or \(50\) displaying the whole range of values for \(C_{d,q,r}\). Subsequently, in Section \ref{subsec:q2,5,10Tails}, we specifically examine the tails of these values for \(q\) in a couple of different dimensions \(d\). Following this, we briefly investigate the special cases of \(q = 1\) (Section~\ref{subsec:q1}) and \(q = 0.5\) (Section~\ref{subsec:q=0.5}).

\subsection{Simulations for \texorpdfstring{$q = 2, 5, 10$}{q 2, 5, 10} and \texorpdfstring{$d = 5, 10, 25, 50$}{d 5,10,25,50}}\label{subsec:q2,5,10,d5,10,25,50}
\begin{figure}[ht]
    \centering
    \begin{minipage}{.45\textwidth}
        \centering
        \includegraphics[width=0.9\linewidth]{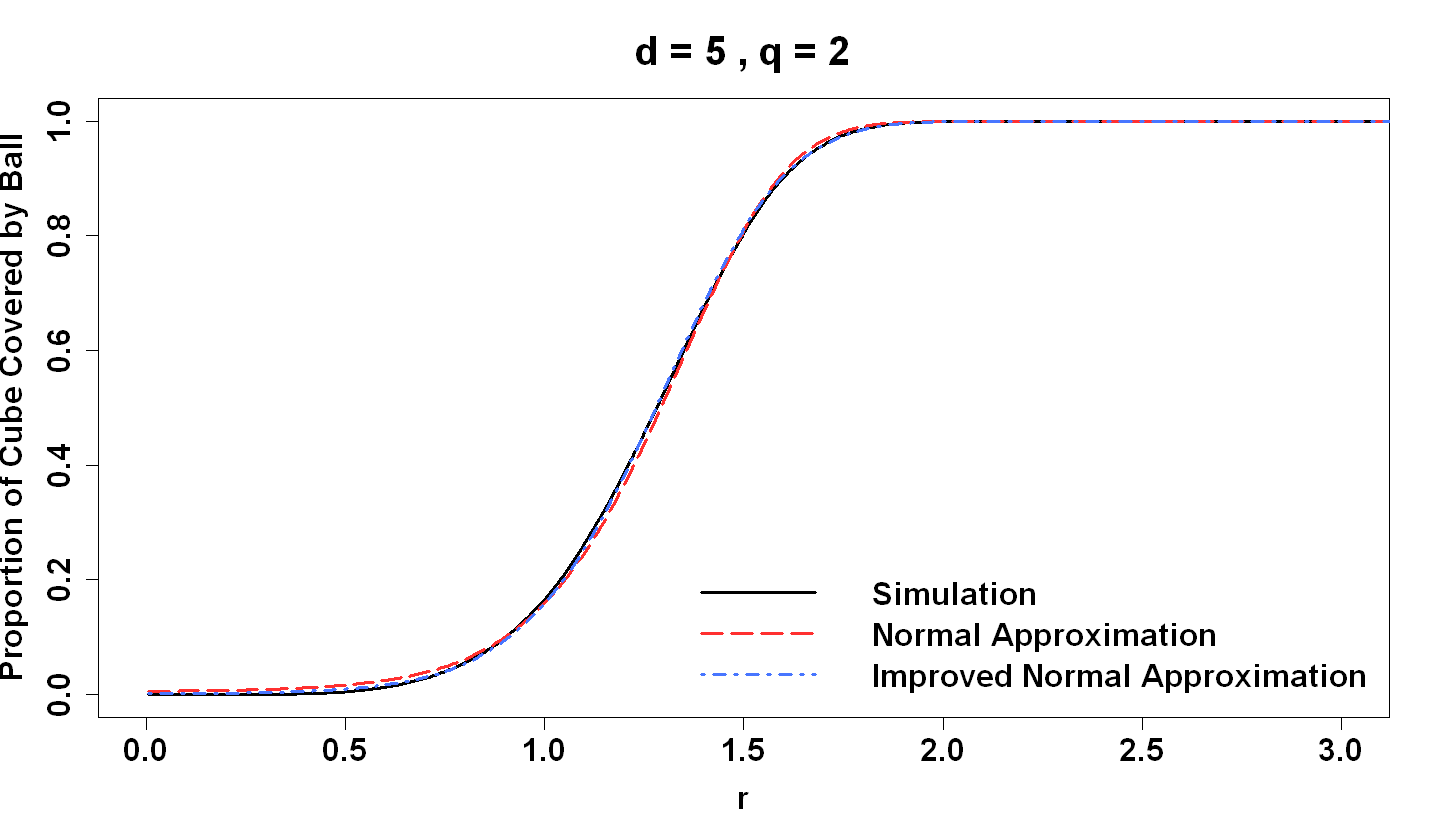}
        \captionsetup{justification=centering}
        \captionof{figure}{\textsl{Graphical representation of Monte Carlo simulation (black) for \(C_{d,q,r}\), \(q=2\), \(d=5\), with Normal (red) and improved Normal (blue) approximations.}}
        \label{fig:q2d5}
        \includegraphics[width=0.9\linewidth]{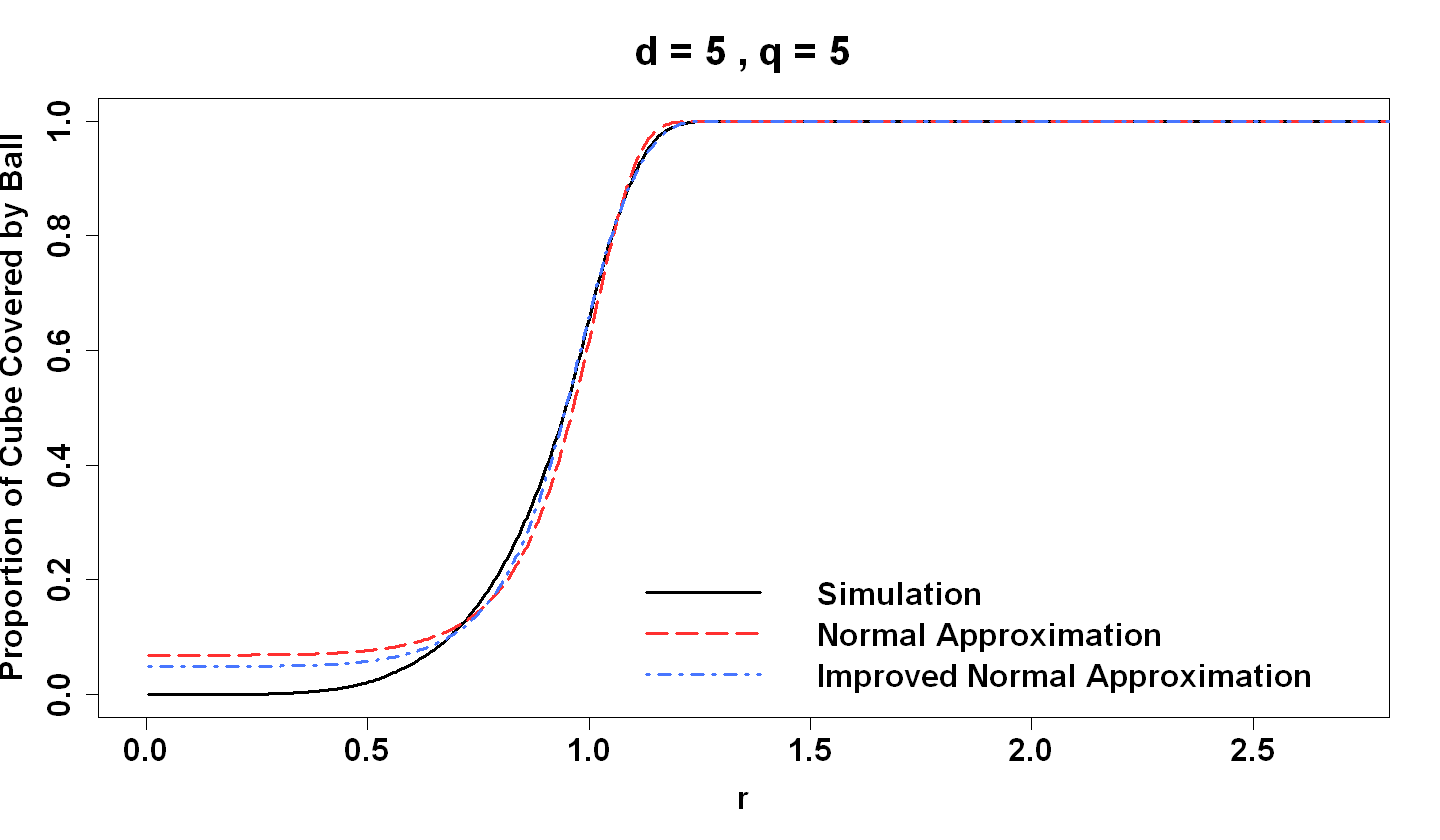}
        \captionsetup{justification=centering}
        \captionof{figure}{\textsl{Graphical representation of Monte Carlo simulation (black) for \(C_{d,q,r}\), \(q=5\), \(d=5\), with Normal (red) and improved Normal (blue) approximations.}}
        \label{fig:q5d5}
    \end{minipage}
    \begin{minipage}{.45\textwidth}
        \centering
        \includegraphics[width=0.9\linewidth]{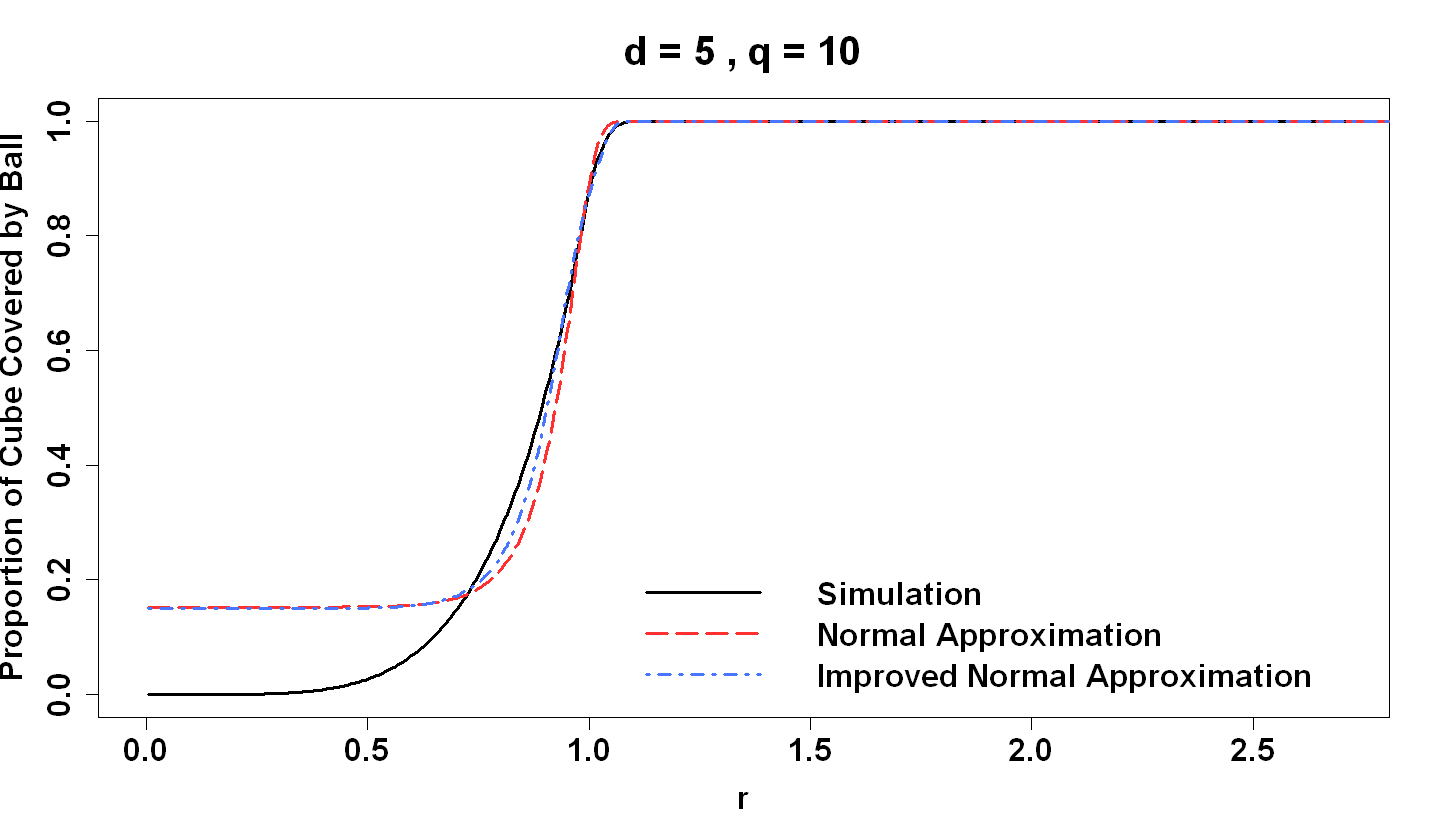}
        \captionsetup{justification=centering}
        \captionof{figure}{\textsl{Graphical representation of Monte Carlo simulation (black) for \(C_{d,q,r}\), \(q=10\), \(d=5\), with Normal (red) and improved Normal (blue) approximations.}}
        \label{fig:q10d5}
        \includegraphics[width=0.9\linewidth]{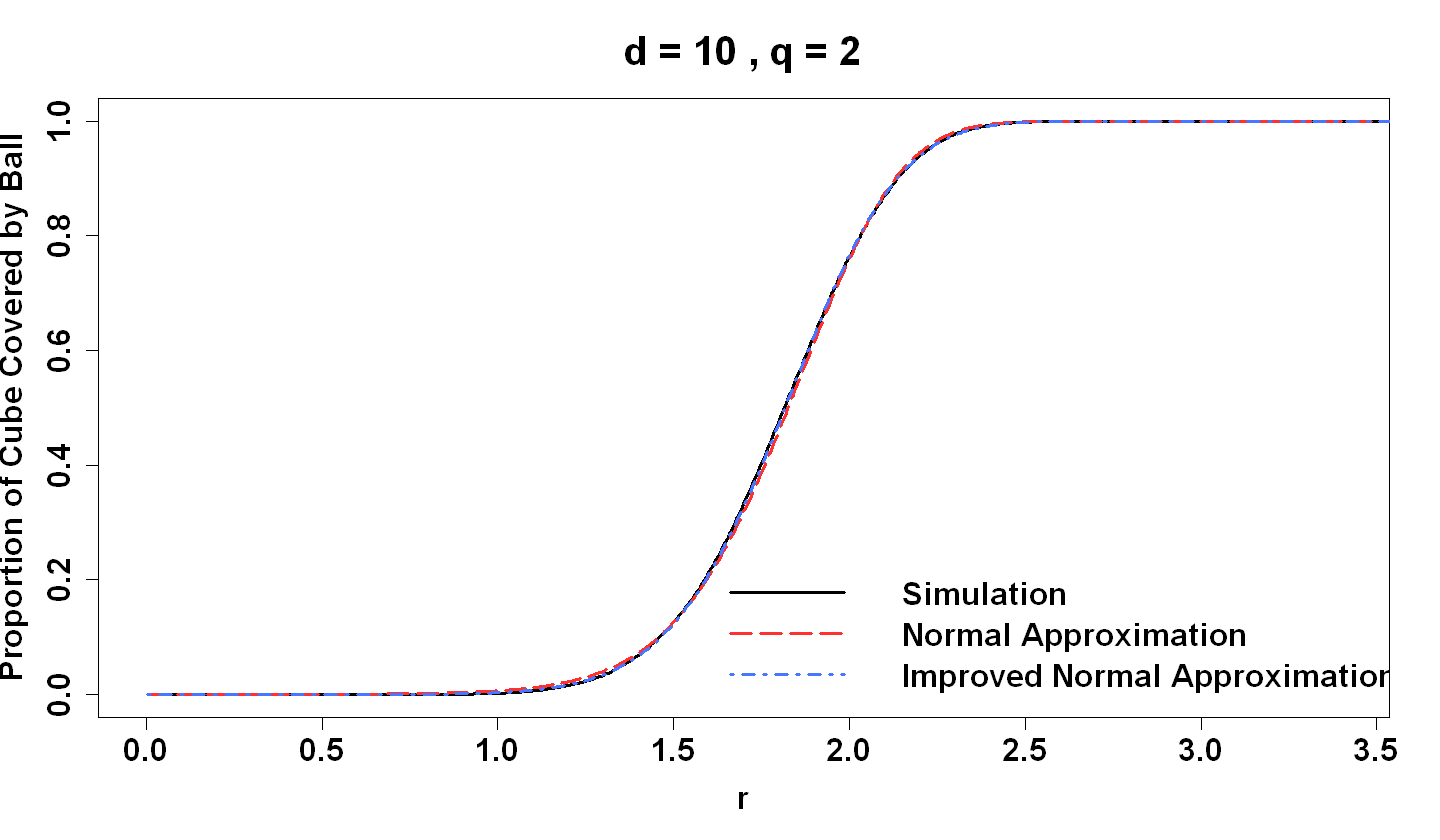}
        \captionsetup{justification=centering}
        \captionof{figure}{\textsl{Graphical representation of Monte Carlo simulation (black) for \(C_{d,q,r}\), \(q=2\), \(d=10\), with Normal (red) and improved Normal (blue) approximations.}}
        \label{fig:q2d10}
    \end{minipage}
\end{figure}

\begin{figure}[ht!]
    \centering
    \begin{minipage}{.45\textwidth}
        \centering
        \includegraphics[width=0.9\linewidth]{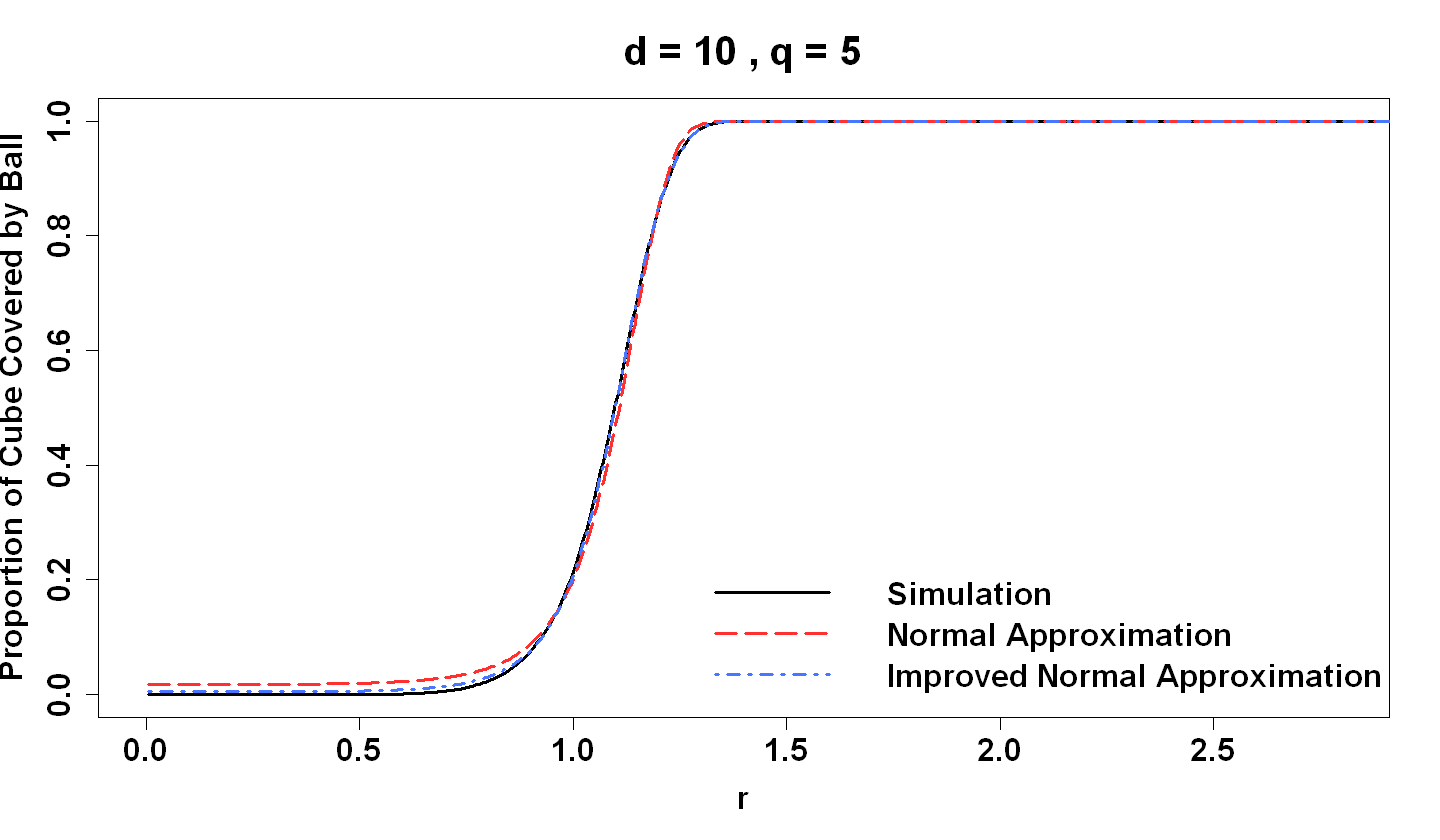}
        \captionsetup{justification=centering}
        \captionof{figure}{\textsl{Graphical representation of Monte Carlo simulation (black) for \(C_{d,q,r}\), \(q=5\), \(d=10\), with Normal (red) and improved Normal (blue) approximations.}}
        \label{fig:q5d10}
        \includegraphics[width=0.9\linewidth]{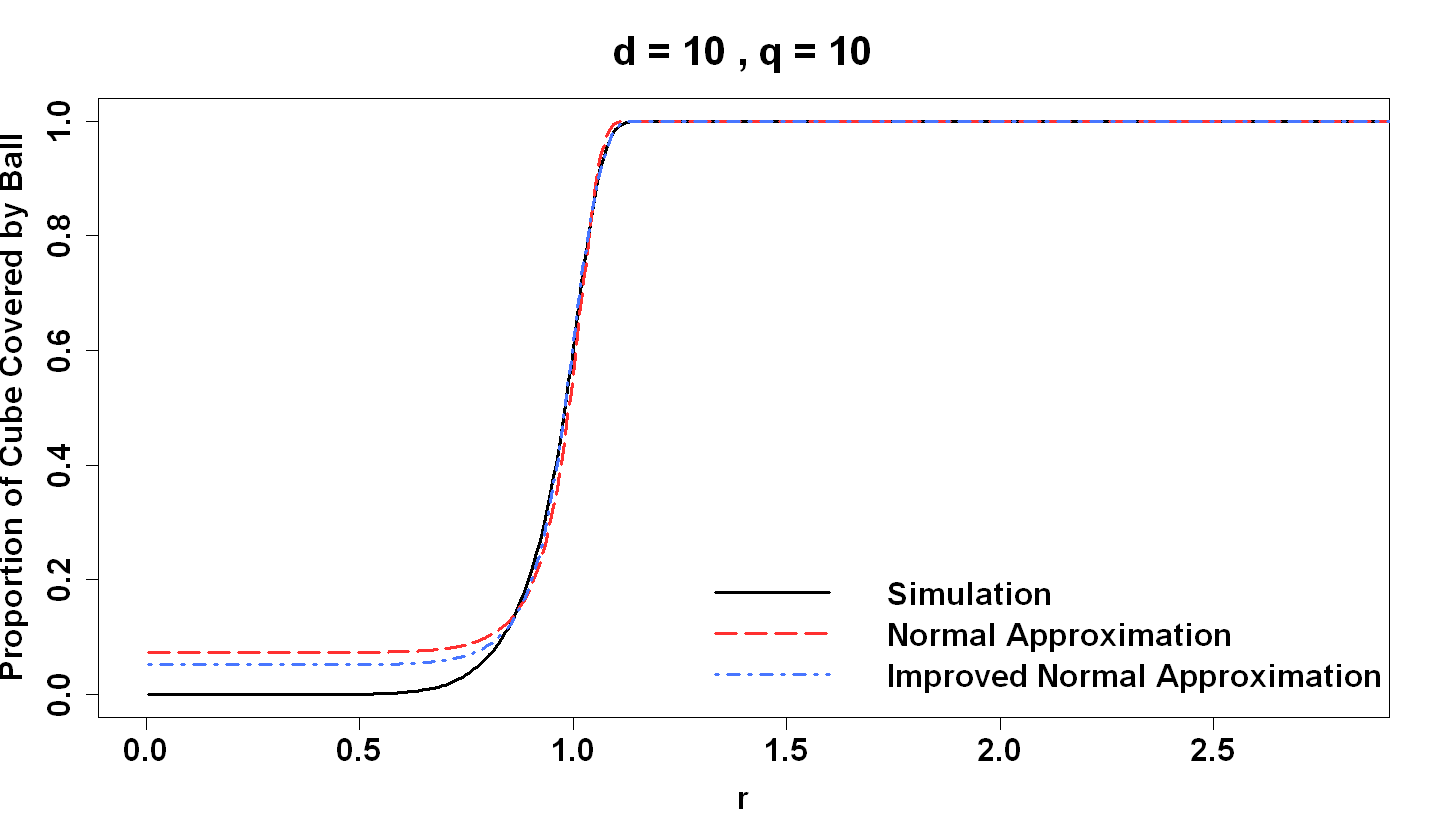}
        \captionsetup{justification=centering}
        \captionof{figure}{\textsl{Graphical representation of Monte Carlo simulation (black) for \(C_{d,q,r}\), \(q=10\), \(d=10\), with Normal (red) and improved Normal (blue) approximations.}}
        \label{fig:q10d10}
        \includegraphics[width=0.9\linewidth]{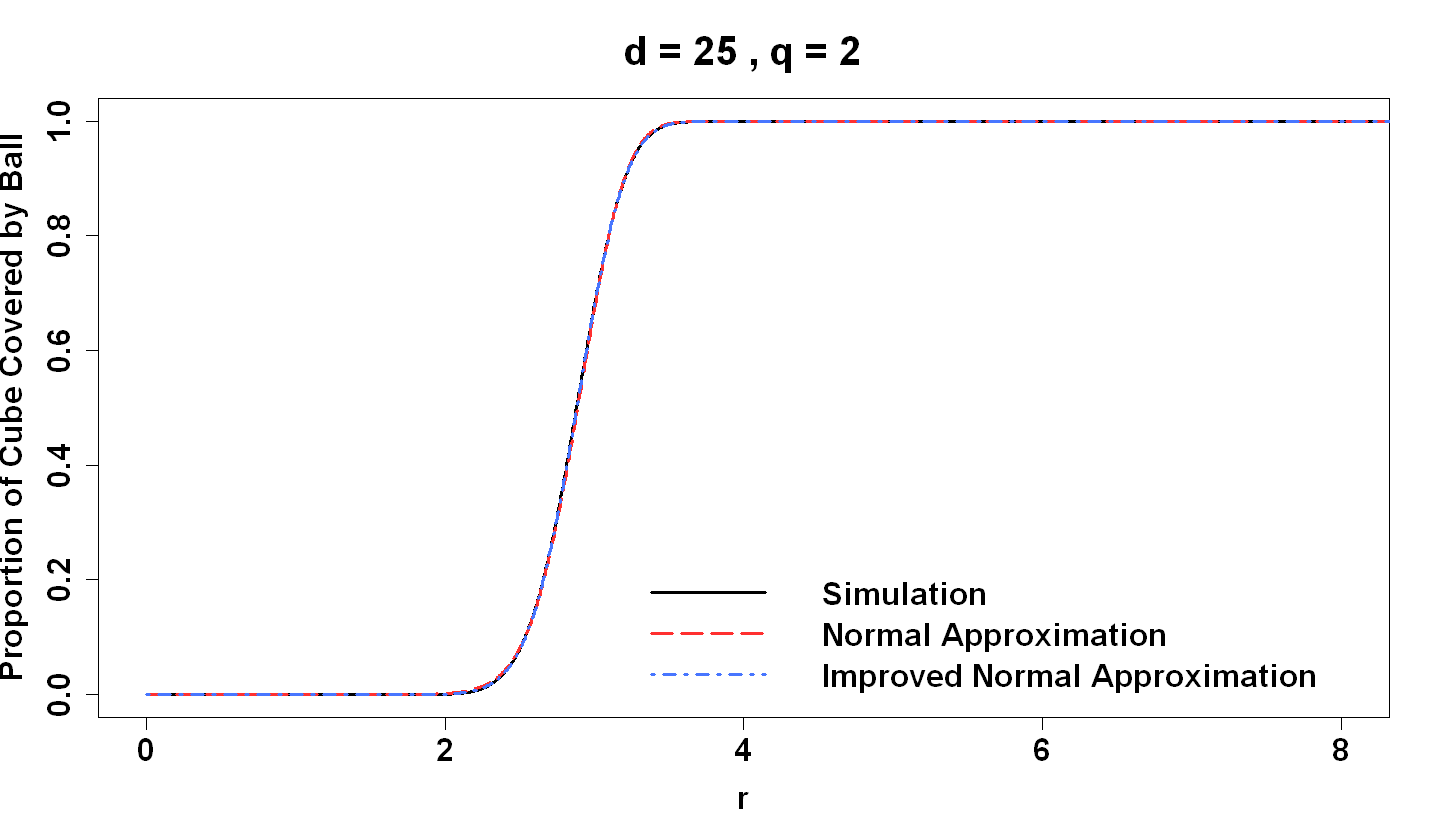}
        \captionsetup{justification=centering}
        \captionof{figure}{\textsl{Graphical representation of Monte Carlo simulation (black) for \(C_{d,q,r}\), \(q=2\), \(d=25\), with Normal (red) and improved Normal (blue) approximations.}}
        \label{fig:q2d25}
        \includegraphics[width=0.9\linewidth]{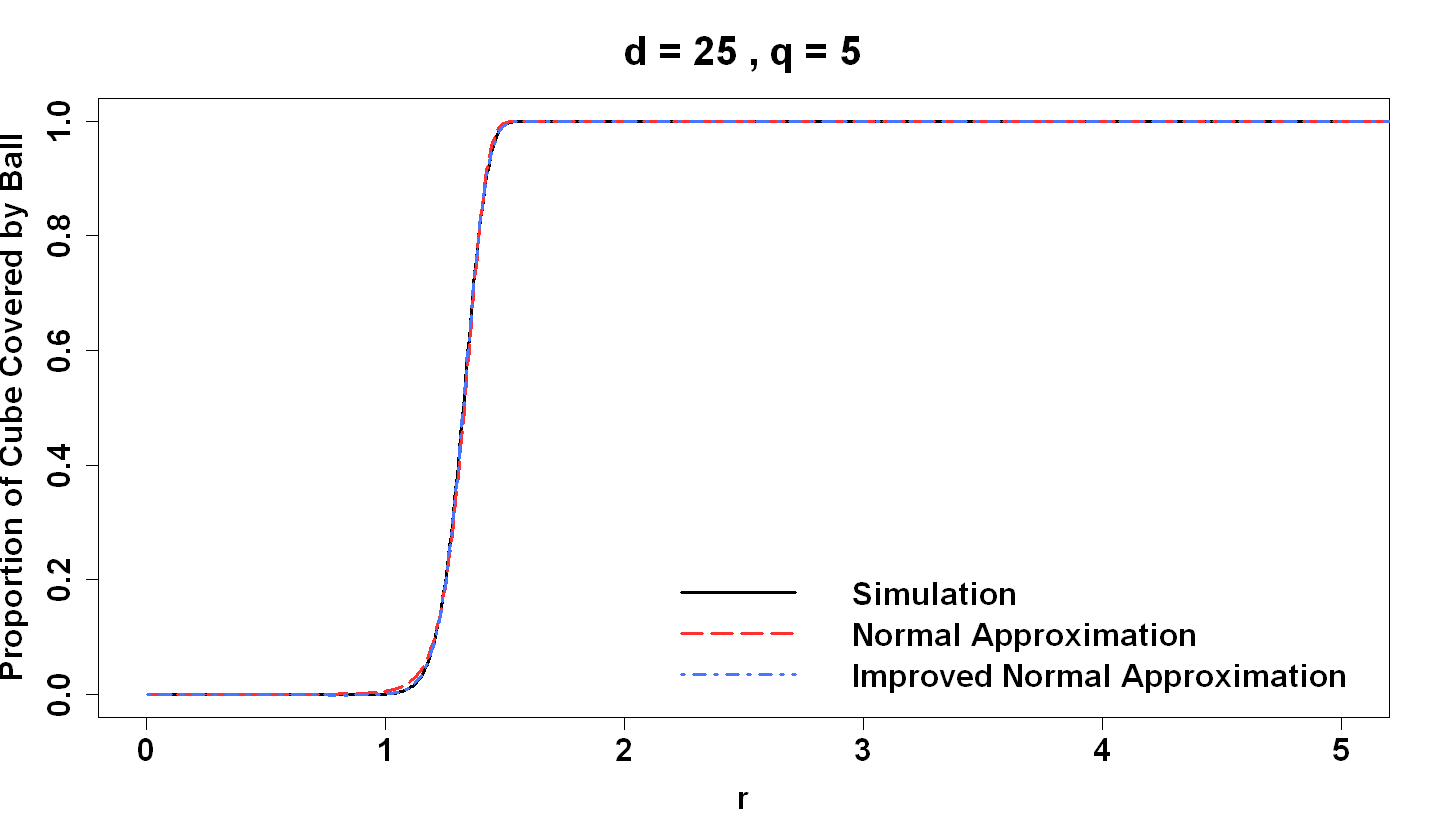}
        \captionsetup{justification=centering}
        \captionof{figure}{\textsl{Graphical representation of Monte Carlo simulation (black) for \(C_{d,q,r}\), \(q=5\), \(d=25\), with Normal (red) and improved Normal (blue) approximations.}}
        \label{fig:q5d25}
    \end{minipage}
    \begin{minipage}{.45\textwidth}
        \centering
        \includegraphics[width=0.9\linewidth]{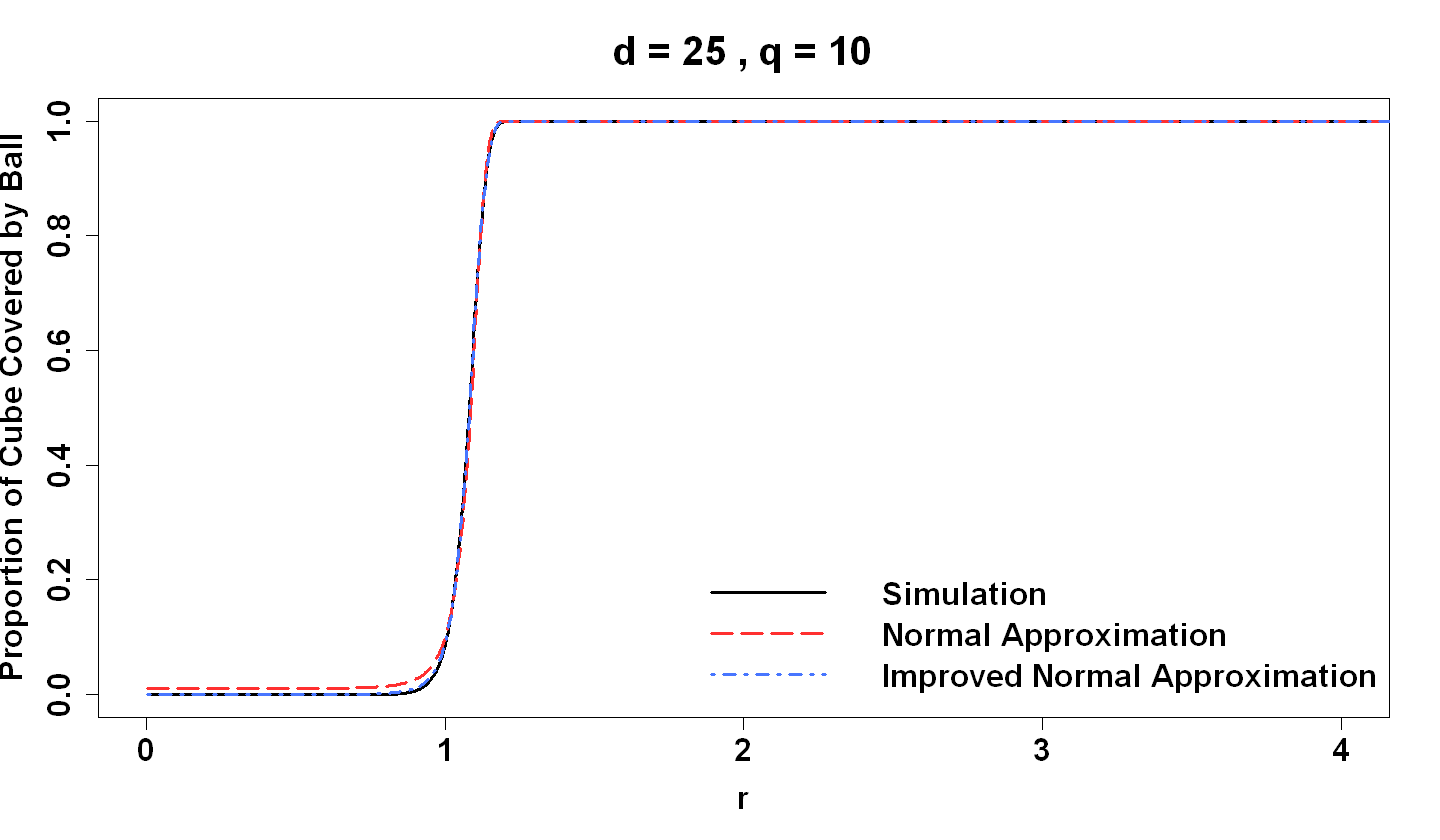}
        \captionsetup{justification=centering}
        \captionof{figure}{\textsl{Graphical representation of Monte Carlo simulation (black) for \(C_{d,q,r}\), \(q=10\), \(d=25\), with Normal (red) and improved Normal (blue) approximations.}}
        \label{fig:q10d25}
        \includegraphics[width=0.9\linewidth]{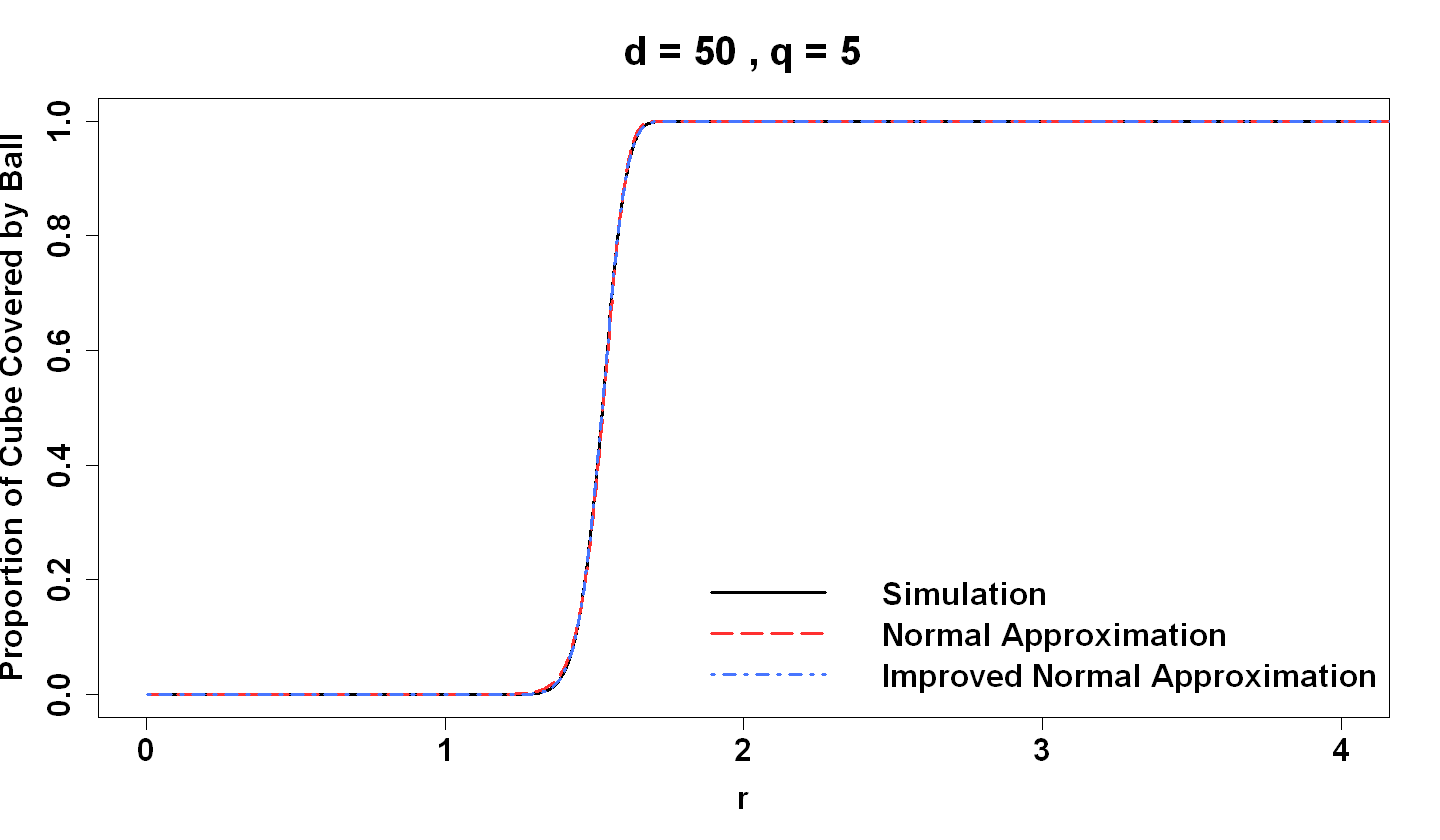}
        \captionsetup{justification=centering}
        \captionof{figure}{\textsl{Graphical representation of Monte Carlo simulation (black) for \(C_{d,q,r}\), \(q=5\), \(d=50\), with Normal (red) and improved Normal (blue) approximations.}}
        \label{fig:q5d50}
        \includegraphics[width=0.9\linewidth]{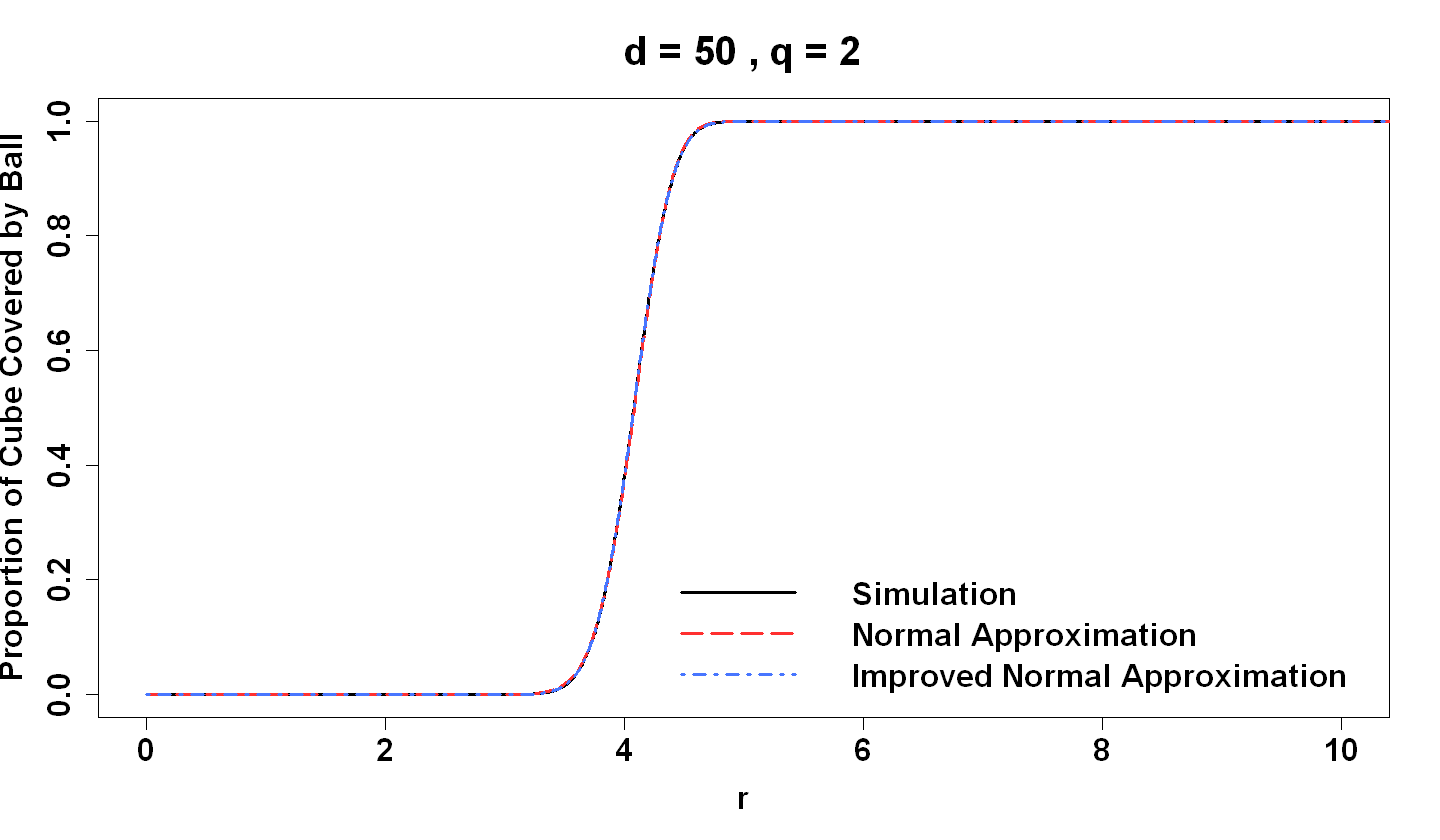}
        \captionsetup{justification=centering}
        \captionof{figure}{\textsl{Graphical representation of Monte Carlo simulation (black) for \(C_{d,q,r}\), \(q=2\), \(d=50\), with Normal (red) and improved Normal (blue) approximations.}}
        \label{fig:q2d50}
        \includegraphics[width=0.9\linewidth]{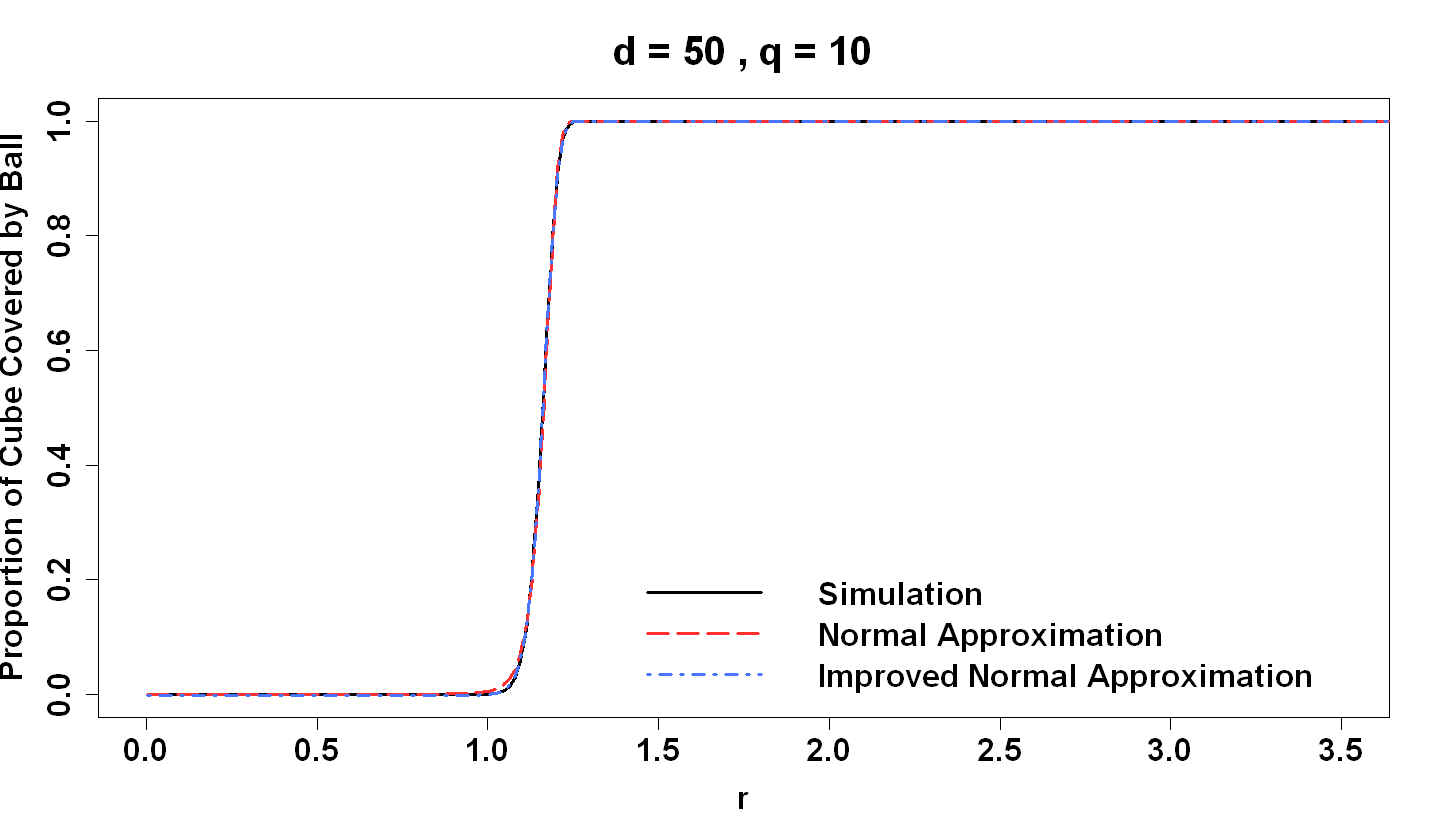}
        \captionsetup{justification=centering}
        \captionof{figure}{\textsl{Graphical representation of Monte Carlo simulation (black) for \(C_{d,q,r}\), \(q=10\), \(d=50\), with Normal (red) and improved Normal (blue) approximations.}}
        \label{fig:q10d50}
    \end{minipage}
\end{figure}

In this section we explore the effect of increasing \(d\) for fixed \(q > 1\), and also of increasing \(q\) when keeping \(d\) constant. The results of this can be seen in Figures \ref{fig:q2d5} - \ref{fig:q10d50}.

As might be expected due to the assumption of high dimension made when deriving the approximations, if we compare graphs vertically for fixed \(q\) we see that the accuracy of both approximations increases as the dimension \(d\) increases for all cases. However, comparing between different values of \(q\) keeping \(d\) the same, we observe that the validity of the approximations decreases as \(q\) increases.

For example, in dimension \(d = 10\), Figure \ref{fig:q2d10} with \(q = 2\) has good agreement between both approximations and the Monte Carlo simulation. Then, in Figure \ref{fig:q5d10} that has \(q = 5\), the improved Normal approximation matches the simulations well whilst the first Normal approximation is less reliable. However, in Figure \ref{fig:q10d10} with \(q = 10\), neither approximation formula fits the curve of the simulation, particularly for small values of \(r\) (e.g. \(r < 1\)).

From Figures \ref{fig:q2d5} - \ref{fig:q10d50}, we may deduce that for the improved Normal approximation to appear very accurate at this view over the whole range, it would seem that we require only a dimension of \(d = 10\) for \(q = 2\). When \(q = 5\), we need to take \(d\) at least \(25\), whilst for \(q = 10\), \(d = 50\) is necessary to obtain a very strong likeness between the curves. Additionally, for the first Normal approximation, we seem to need a dimension of at least \(25\) for \(q = 2\), \(50\) for \(q = 5\) and some value greater than \(50\) for \(q = 10\) to have high proximity to the simulation over the whole range.

Using these insights, we may conclude that the accuracy of both approximation formulae increases as the dimension gets larger. However, higher values of \(d\) are required for the first approximation to appear valid when compared to the improved Normal approximation. Additionally, it appears that the higher the value of \(q\), the larger the dimension needed to ensure the same level of accuracy.


\cleardoublepage
\subsection{Simulations for the tails of \texorpdfstring{$q = 2, 5, 10$}{q 2, 5, 10}}\label{subsec:q2,5,10Tails}
Having explored the general trend for \(q = 2, 5, 10\) and \(d = 5, 10, 25, 50\) over the whole range of values for \(C_{d,q,r}\) in Section \ref{subsec:q2,5,10,d5,10,25,50}, we now examine the tails of two distributions for  \(q = 2, 5, 10\) in more detail.

\begin{figure}[ht!]
    \centering
    \begin{minipage}{.45\textwidth}
        \centering
        \includegraphics[width=0.9\linewidth]{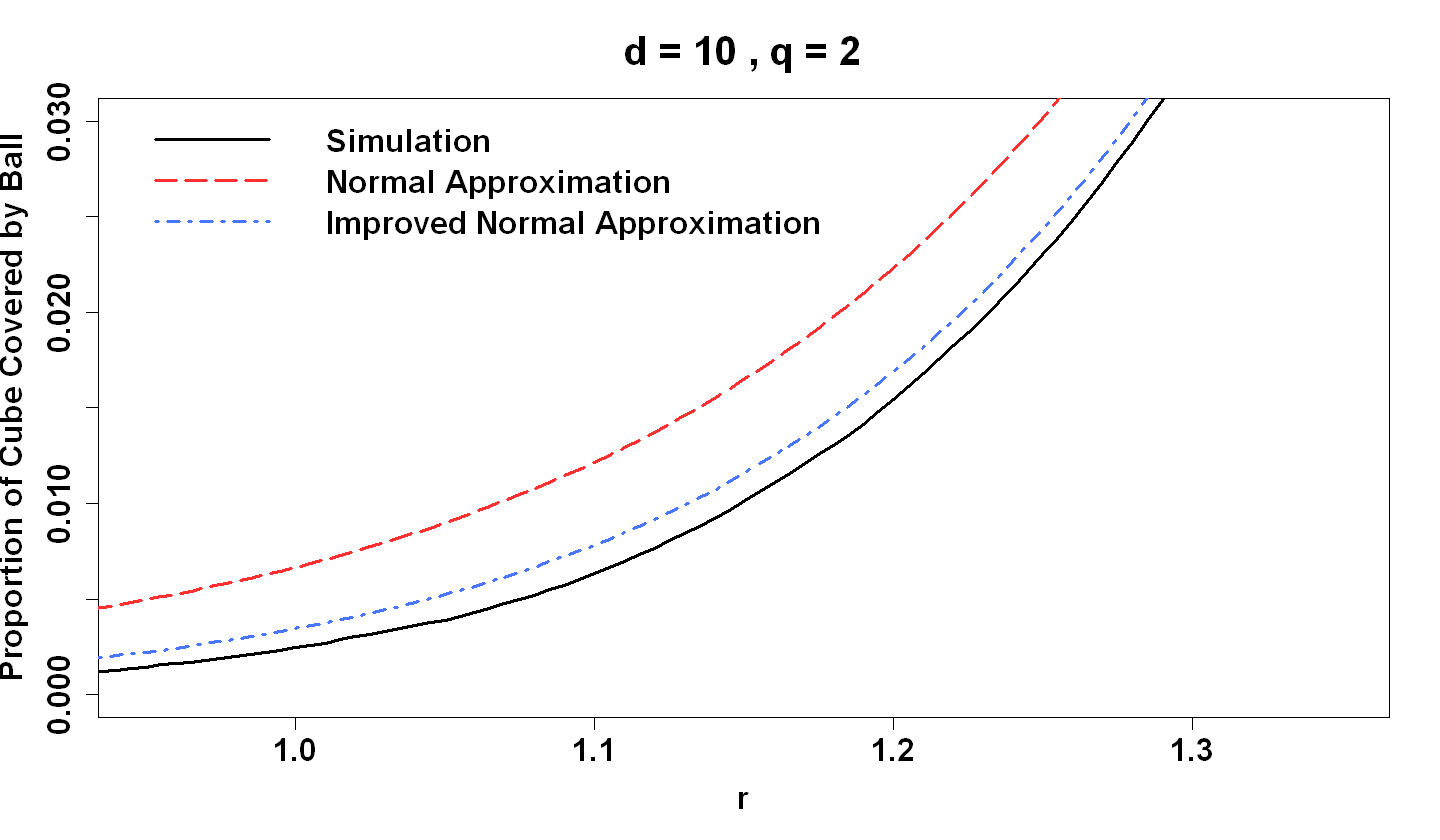}
        \captionsetup{justification=centering}
        \captionof{figure}{\textsl{Graphical representation of Monte Carlo simulation (black) for \(C_{d,q,r}\), \(q=2\), \(d=10\), with Normal (red) and improved Normal (blue) approximations, zoomed in on the tail.}}
        \label{fig:q2d10Tail}
        \includegraphics[width=0.9\linewidth]{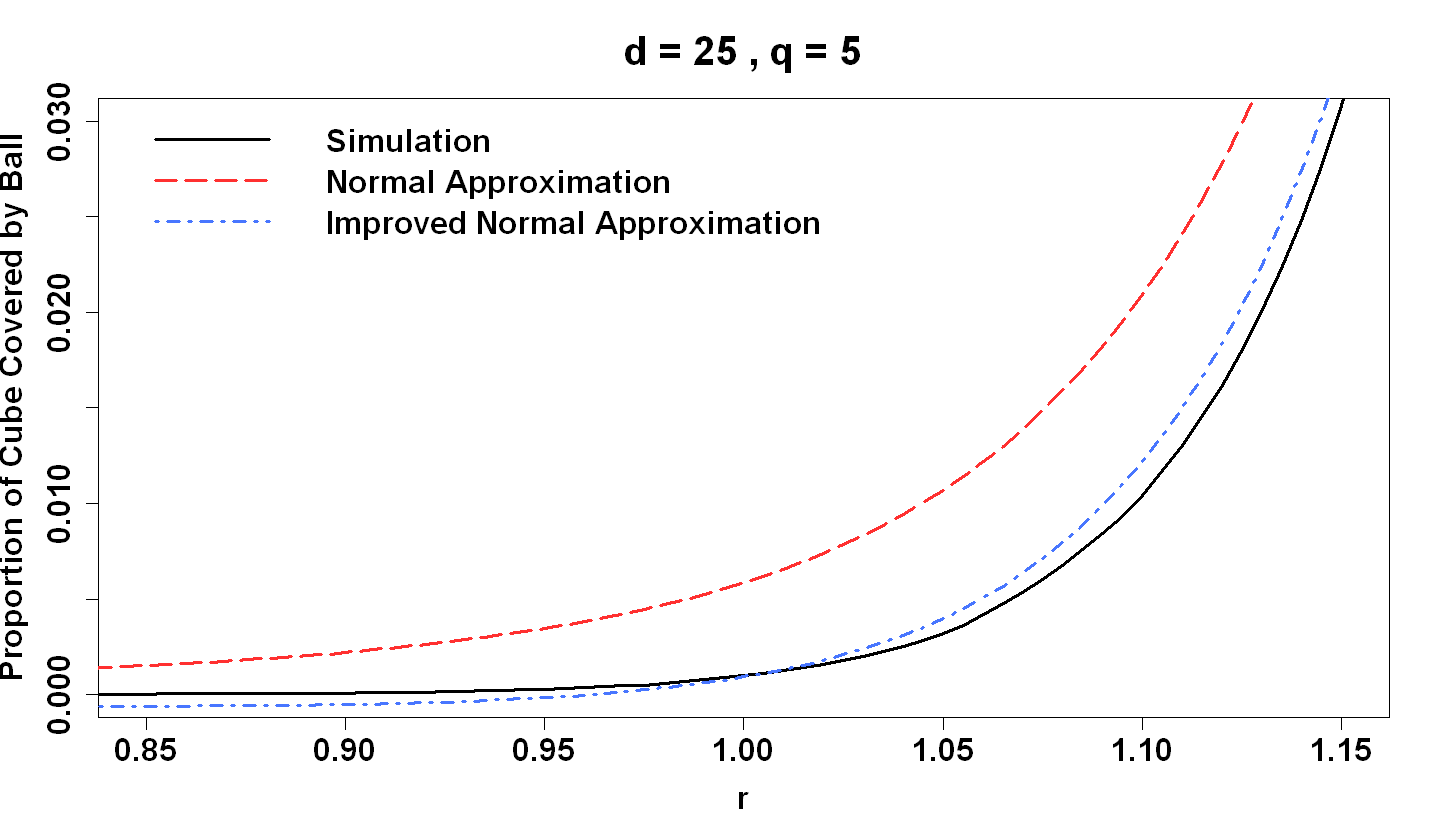}
        \captionsetup{justification=centering}
        \captionof{figure}{\textsl{Graphical representation of Monte Carlo simulation (black) for \(C_{d,q,r}\), \(q=5\), \(d=25\), with Normal (red) and improved Normal (blue) approximations, zoomed in on the tail.}}
        \label{fig:q5d25Tail}
        \includegraphics[width=0.9\linewidth]{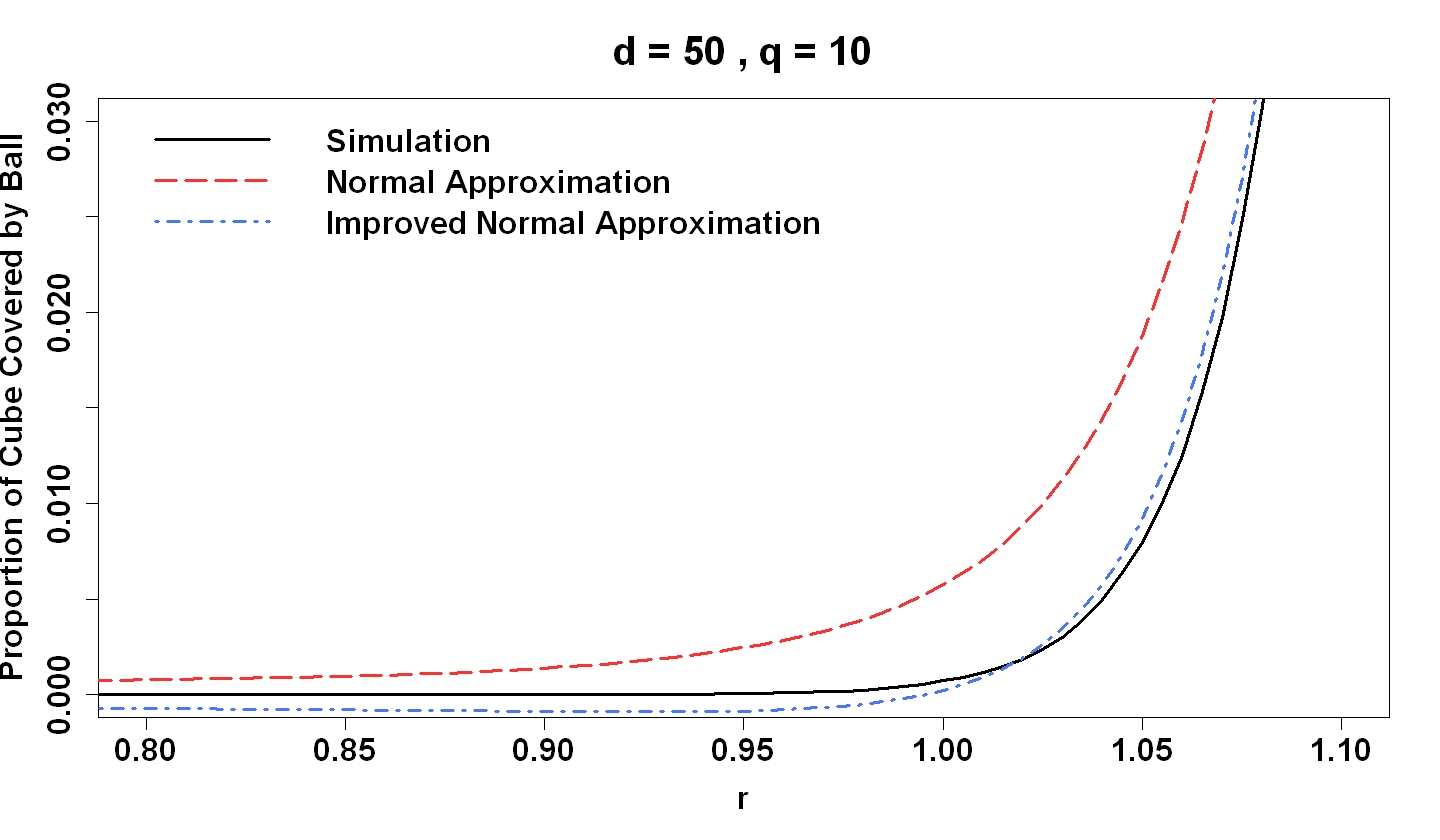}
        \captionsetup{justification=centering}
        \captionof{figure}{\textsl{Graphical representation of Monte Carlo simulation (black) for \(C_{d,q,r}\), \(q=10\), \(d=50\), with Normal (red) and improved Normal (blue) approximations, zoomed in on the tail.}}
        \label{fig:q10d50Tail}
    \end{minipage}
    \begin{minipage}{.45\textwidth}
        \centering
        \includegraphics[width=0.9\linewidth]{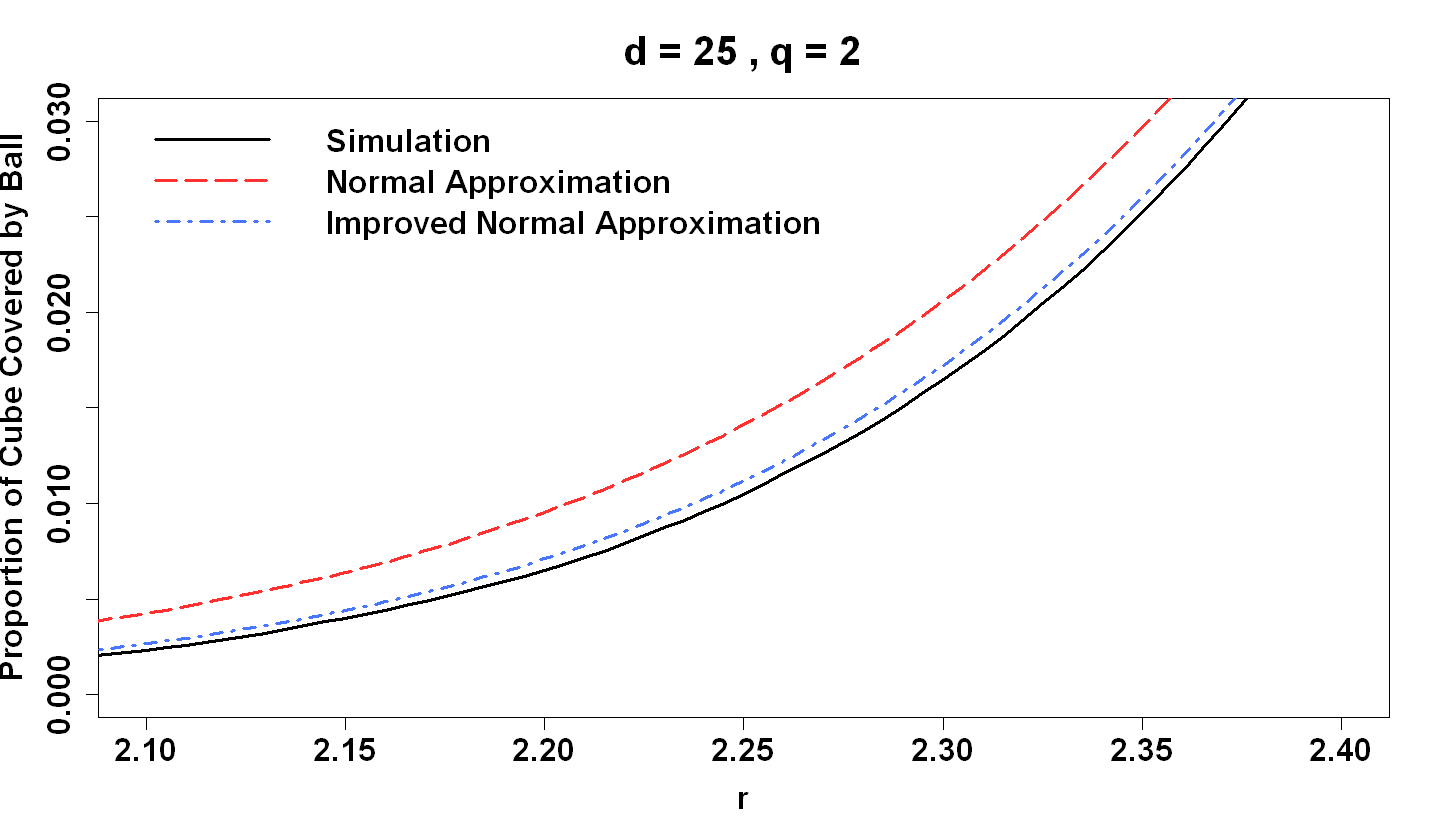}
        \captionsetup{justification=centering}
        \captionof{figure}{\textsl{Graphical representation of Monte Carlo simulation (black) for \(C_{d,q,r}\), \(q=2\), \(d=25\), with Normal (red) and improved Normal (blue) approximations, zoomed in on the tail.}}
        \label{fig:q2d25Tail}
        \includegraphics[width=0.9\linewidth]{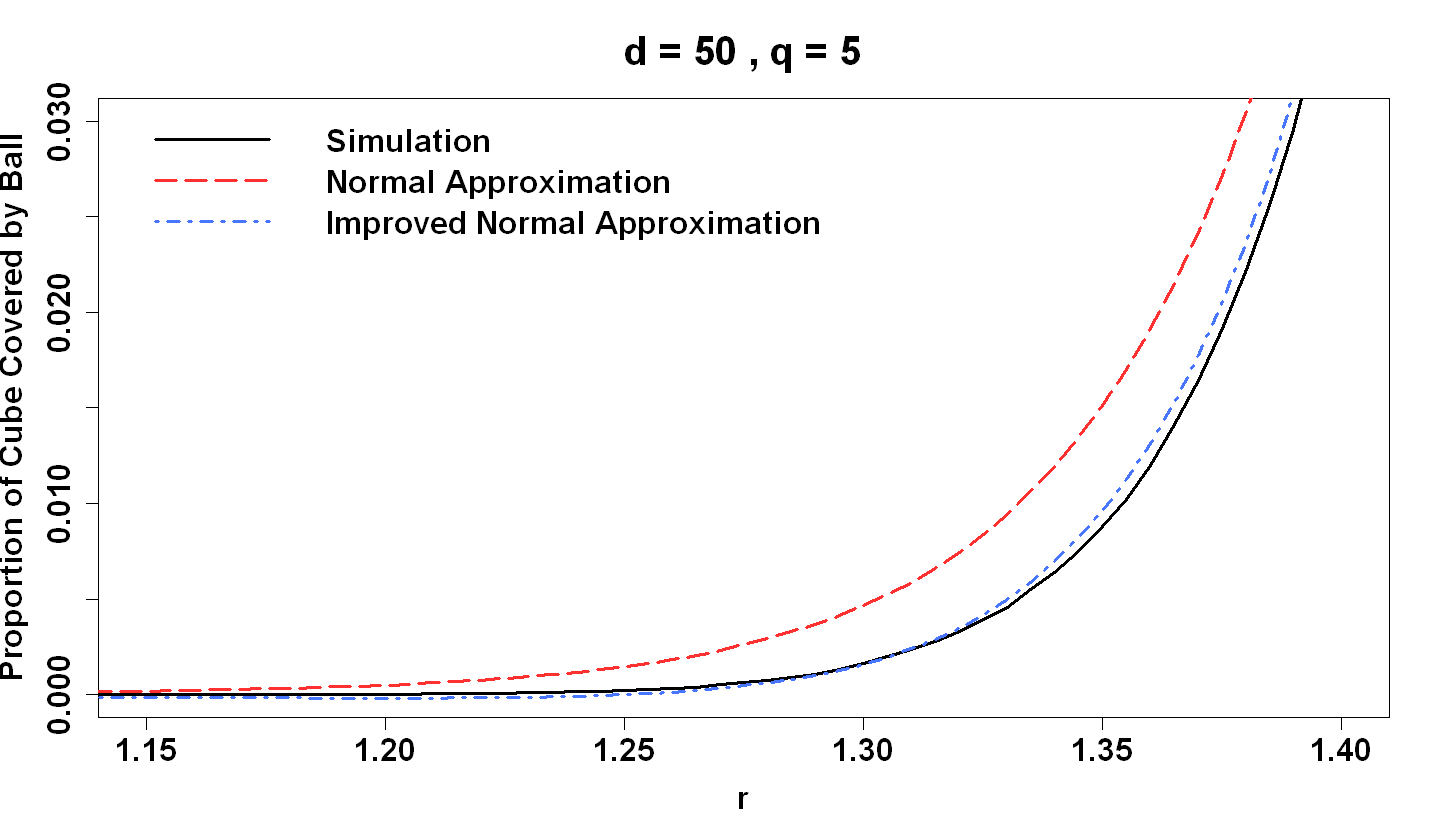}
        \captionsetup{justification=centering}
        \captionof{figure}{\textsl{Graphical representation of Monte Carlo simulation (black) for \(C_{d,q,r}\), \(q=5\), \(d=50\), with Normal (red) and improved Normal (blue) approximations, zoomed in on the tail.}}
        \label{fig:q5d50Tail}
        \includegraphics[width=0.9\linewidth]{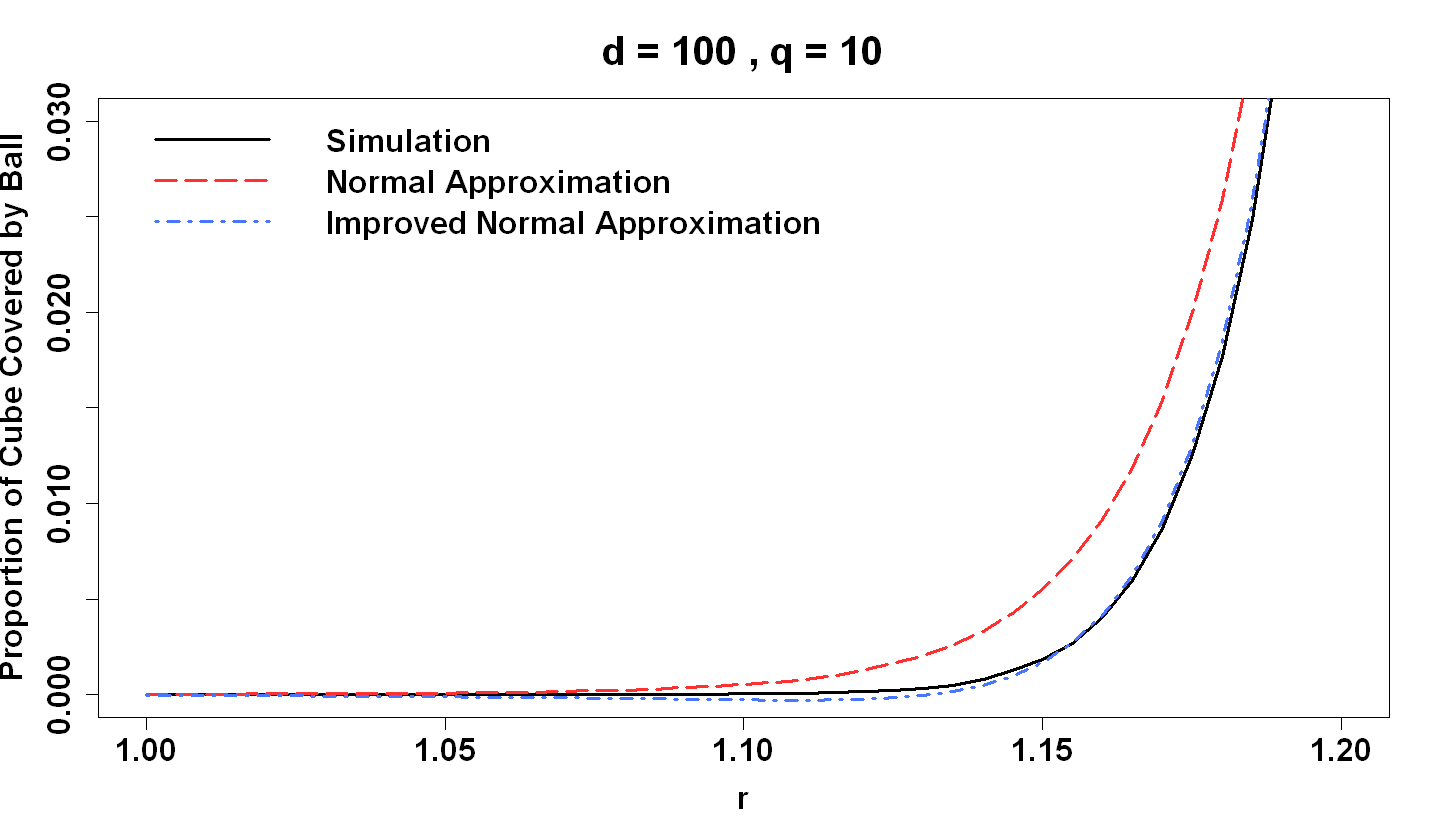}
        \captionsetup{justification=centering}
        \captionof{figure}{\textsl{Graphical representation of Monte Carlo simulation (black) for \(C_{d,q,r}\), \(q=10\), \(d=100\), with Normal (red) and improved Normal (blue) approximations, zoomed in on the tail.}}
        \label{fig:q10d100Tail}
    \end{minipage}
\end{figure}

Contrary to Figure \ref{fig:q2d10} in Section \ref{subsec:q2,5,10,d5,10,25,50}, where the three lines appear almost exactly aligned, we can observe in Figure \ref{fig:q2d10Tail} that there is a difference between the approximations and the Monte Carlo simulation for \(q = 2\) and \(d = 10\), particularly for the first approximation. However, we note that the scale on the \(y\)-axis is small, confirming that the difference would be barely noticeable on a larger scale. Comparing Figure \ref{fig:q2d10Tail} to Figure \ref{fig:q2d25Tail}, where we have increased \(d\) but kept \(q\) the same, we can remark that, in agreement with our suggestions in Section \ref{subsec:q2,5,10,d5,10,25,50}, the gap between the approximation formulae and the simulations has reduced, increasing accuracy of the approximations.

On examining Figures \ref{fig:q5d25Tail}, \ref{fig:q10d50Tail}, \ref{fig:q5d50Tail} and \ref{fig:q10d100Tail}, we may notice an inaccuracy with the improved Normal approximation formula, in that it can give negative values for \(C_{d,q,r}\) which should be a non-negative quantity by the theory presented in Section \ref{sec:Theory}. This negativity is most visible in Figure \ref{fig:q10d100Tail}, however, we again comment that, due to the scale of these graphs, the negative values appear to be close to zero. Since \(d\) has been chosen large in all these cases, despite the negative values, there is strong agreement between the improved Normal approximation formula and the Monte Carlo simulation. To explore this issue of negative values further, it could be investigated whether the second approximation formula could be made more accurate by, for example, taking the absolute value or setting negative values to zero.

\subsection{Simulations for \texorpdfstring{$q = 1$}{q 1} and \texorpdfstring{$d = 10, 25$}{d 10, 25}}\label{subsec:q1}
\begin{figure}[ht]
    \centering
    \begin{minipage}{.45\textwidth}
        \centering
        \includegraphics[width=0.9\linewidth]{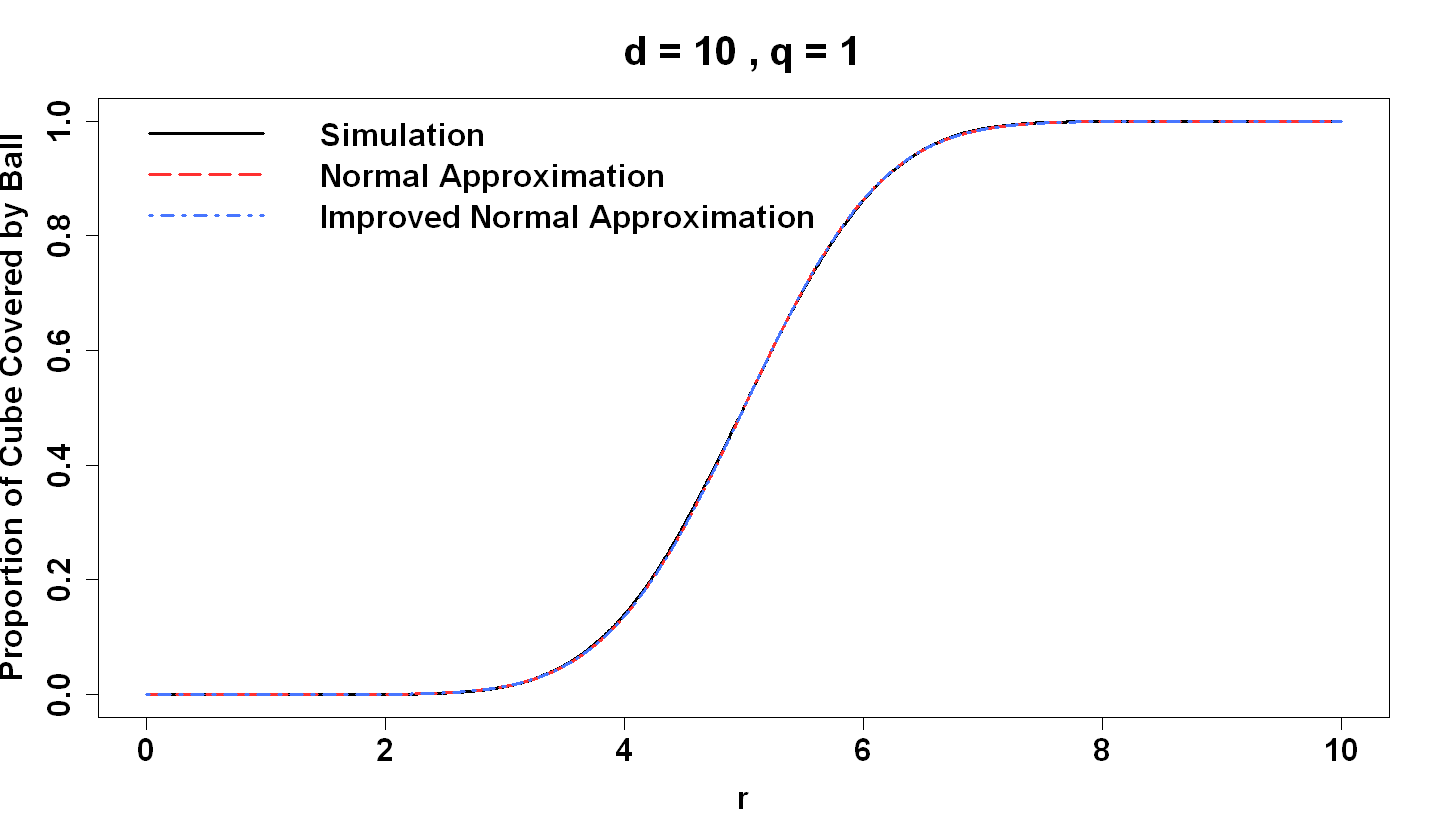}
        \captionsetup{justification=centering}
        \captionof{figure}{\textsl{Graphical representation of Monte Carlo simulation (black) for \(C_{d,q,r}\), \(q=1\), \(d=10\), with Normal (red) and improved Normal (blue) approximations.}}
        \label{fig:q1d10}
        \includegraphics[width=0.9\linewidth]{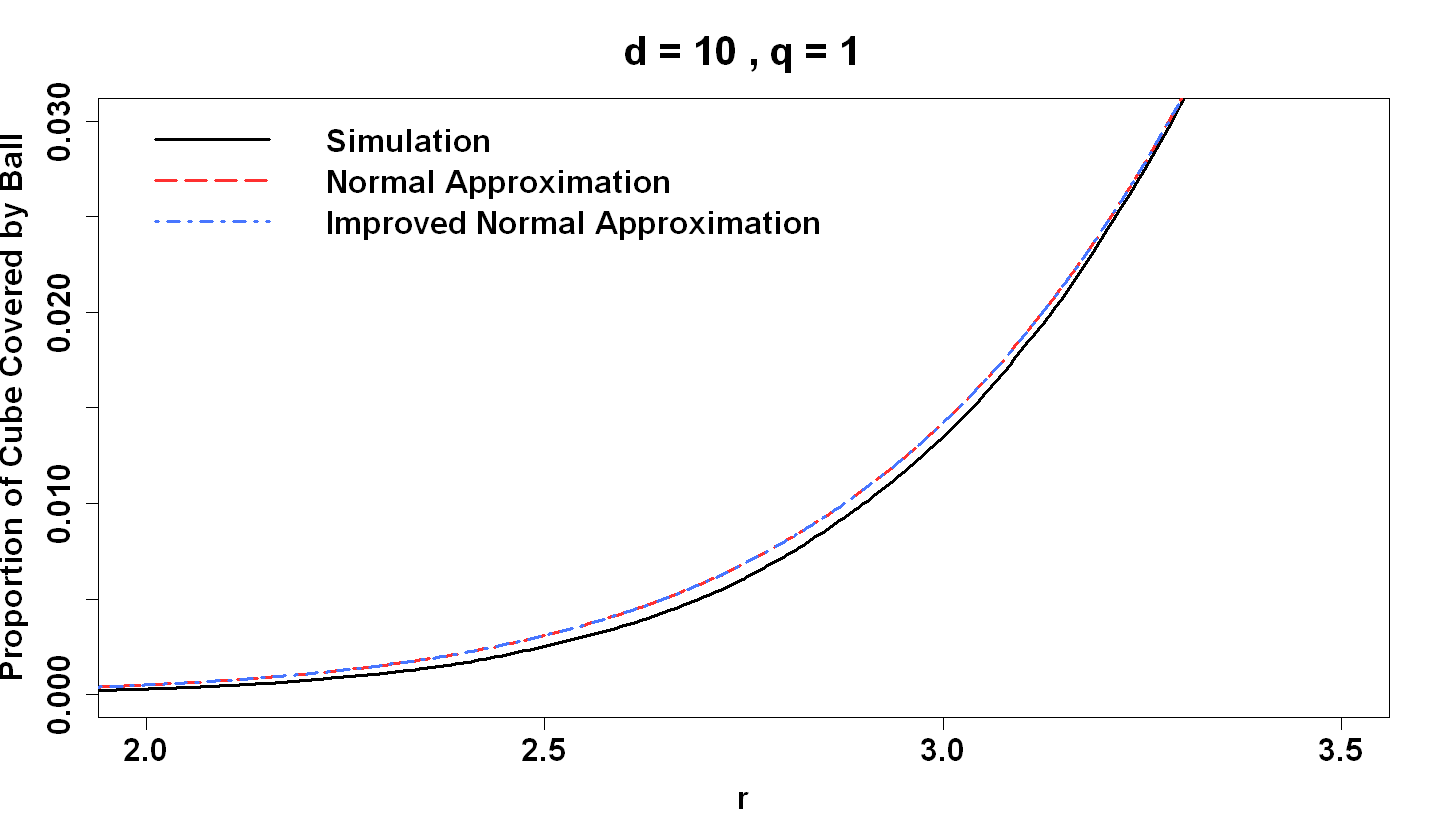}
        \captionsetup{justification=centering}
        \captionof{figure}{\textsl{Graphical representation of Monte Carlo simulation (black) for \(C_{d,q,r}\), \(q=1\), \(d=10\), with Normal (red) and improved Normal (blue) approximations, zoomed in on the tail.}}
        \label{fig:q1d10zoom}
    \end{minipage}
    \begin{minipage}{.45\textwidth}
        \centering
        \includegraphics[width=0.9\linewidth]{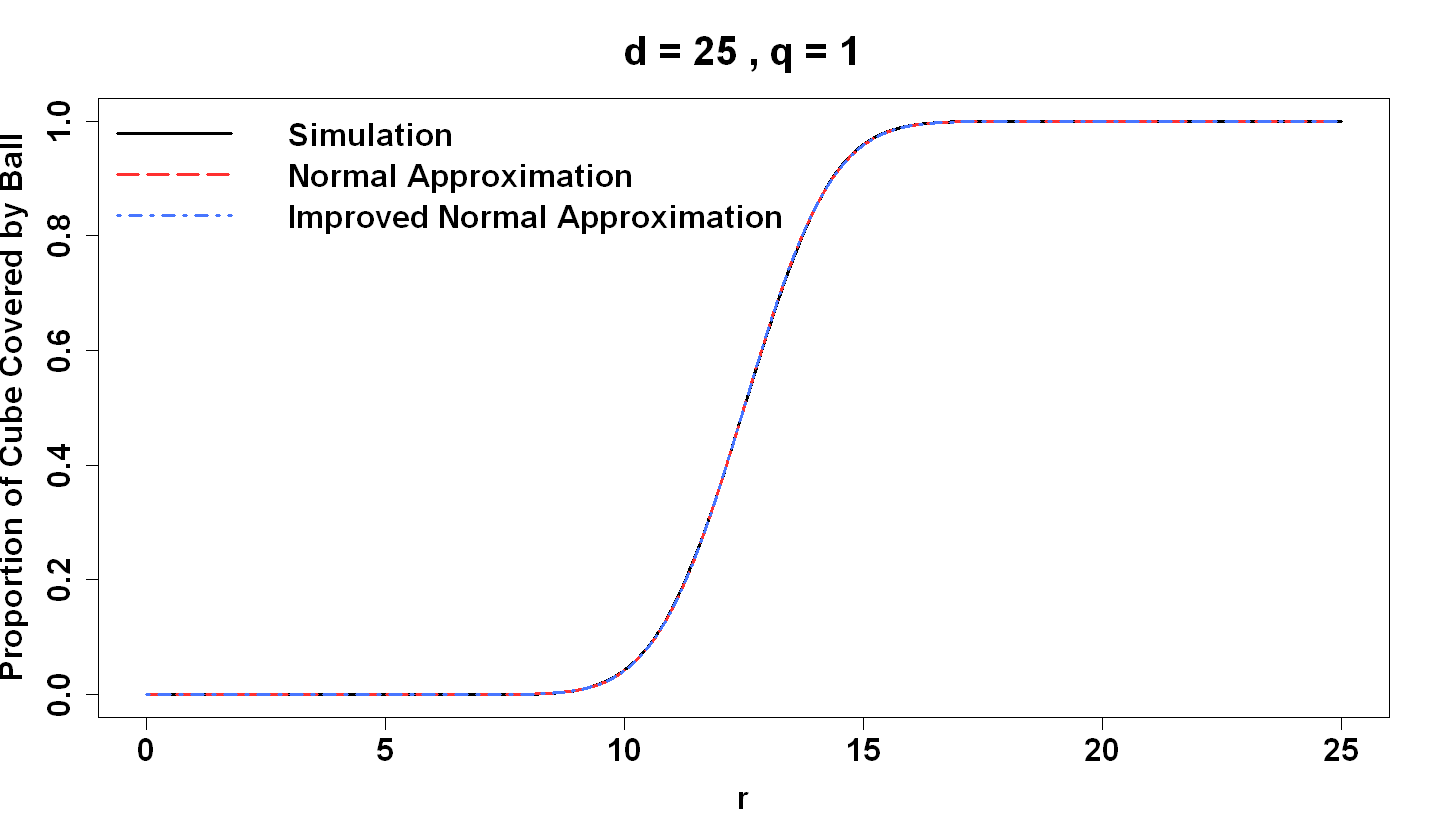}
        \captionsetup{justification=centering}
        \captionof{figure}{\textsl{Graphical representation of Monte Carlo simulation (black) for \(C_{d,q,r}\), \(q=1\), \(d=25\), with Normal (red) and improved Normal (blue) approximations.}}
        \label{fig:q1d25}
        \includegraphics[width=0.9\linewidth]{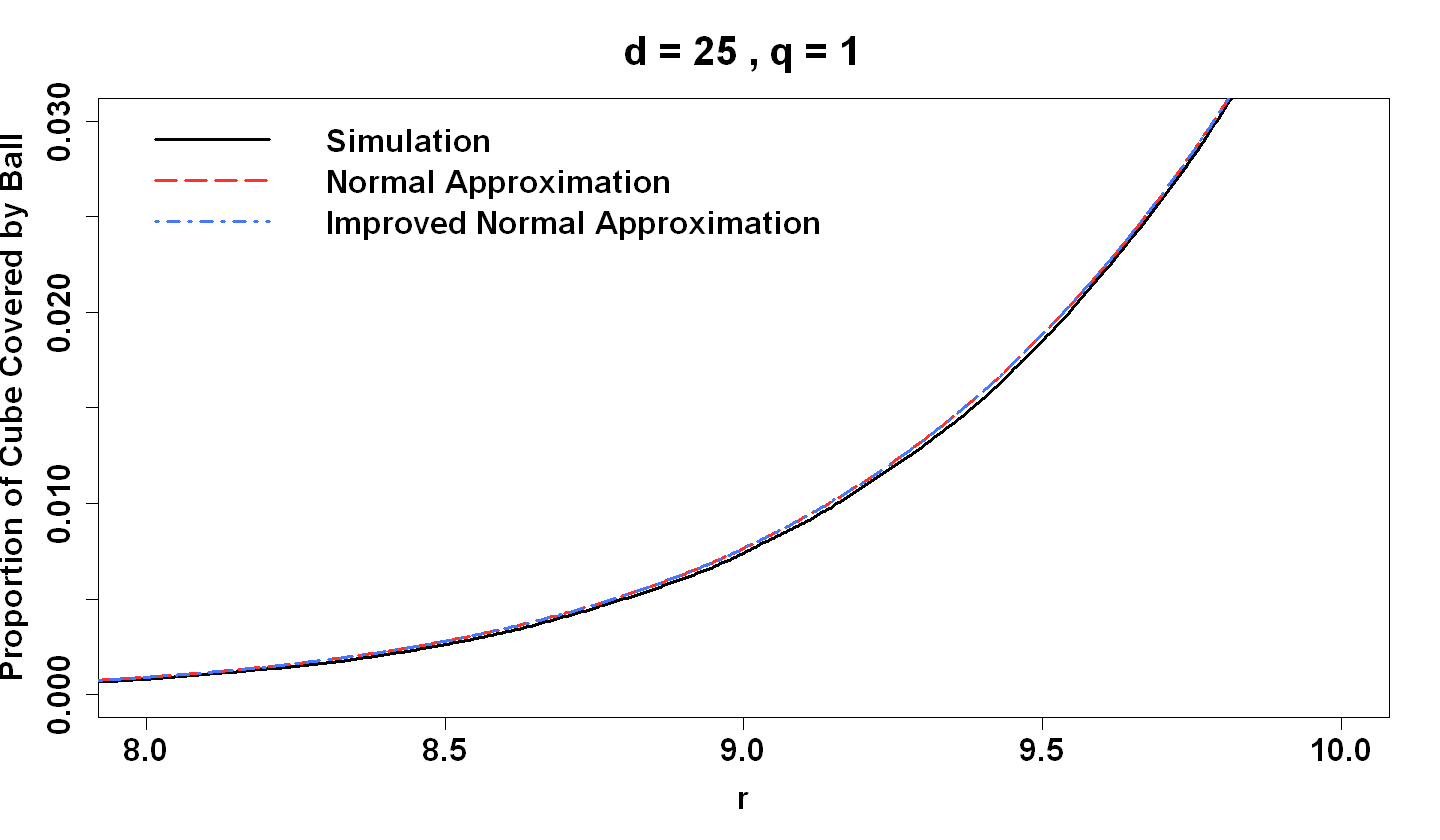}
        \captionsetup{justification=centering}
        \captionof{figure}{\textsl{Graphical representation of Monte Carlo simulation (black) for \(C_{d,q,r}\), \(q=1\), \(d=25\), with Normal (red) and improved Normal (blue) approximations, zoomed in on the tail.}}
        \label{fig:q1d25zoom}
    \end{minipage}
\end{figure}

In this section we consider the case of \(q = 1\) in the dimensions \(d = 10\) and \(d = 25\). By examining the formula in equation (\ref{eq:Approx2BoxVersion}), we can observe that the second term in the improved Normal approximation is zero for this scenario, due to the factor \((q - 1)\). Additionally, by comparing equation (\ref{eq:Approx2BoxVersiontrValue}) with the argument of the standard Normal CDF, \(\Phi\) in equation (\ref{eq:Approx1BoxVersion}), we conclude that when \(q = 1\) both approximation formulae are the same. Accordingly, we see that the red and blue lines exactly coincide in all of Figures \ref{fig:q1d10}, \ref{fig:q1d10zoom}, \ref{fig:q1d25} and \ref{fig:q1d25zoom}.

Studying Figure \ref{fig:q1d10}, we see that over the whole range of \(r\) for \(d = 10\), there seems to be good agreement between the approximations and the simulation. When \(d\) increases to \(25\) in Figure \ref{fig:q1d25}, this is even more apparent. Focusing in on the tail, we see in Figure \ref{fig:q1d10zoom} that there is some difference between the approximations and simulations for \(d = 10\). However, at dimension \(d = 25\) the red and blues lines of the approximations are almost completely aligned with the black line of the Monte Carlo simulation. Hence, we may deduce that \(d = 25\) already appears to be a large enough dimension for the approximations to be a good representation of the true value here. We also suggest that in the case of \(q = 1\), we could just use the first Normal approximation, since the second gives us no gain in accuracy. Further investigations on this topic could explore whether adding the next term in the Edgeworth Expansion for \(q = 1\) would drive any additional improvements in accuracy.

\subsection{Simulations for \texorpdfstring{$q = 0.5$}{q 0.5} and \texorpdfstring{$d = 5, 10$}{d 5, 10}}\label{subsec:q=0.5}
In this final section on simulations we consider \(q=0.5\) as a representation of the case \(0 < q < 1\). Note that in this scenario the \(L_{q}\)-norm is not a true norm \cite{WikiLp}. For example, in dimension \(d=2\), the \(L_{q}\)-ball is no longer a convex object but instead defines a concave body. As a consequence, we may expect that the results we obtain will be different to those from Sections \ref{subsec:q2,5,10,d5,10,25,50} - \ref{subsec:q1}.

\begin{figure}[ht]
\centering
\begin{minipage}{.45\textwidth}
  \centering
  \includegraphics[width=0.9\linewidth]{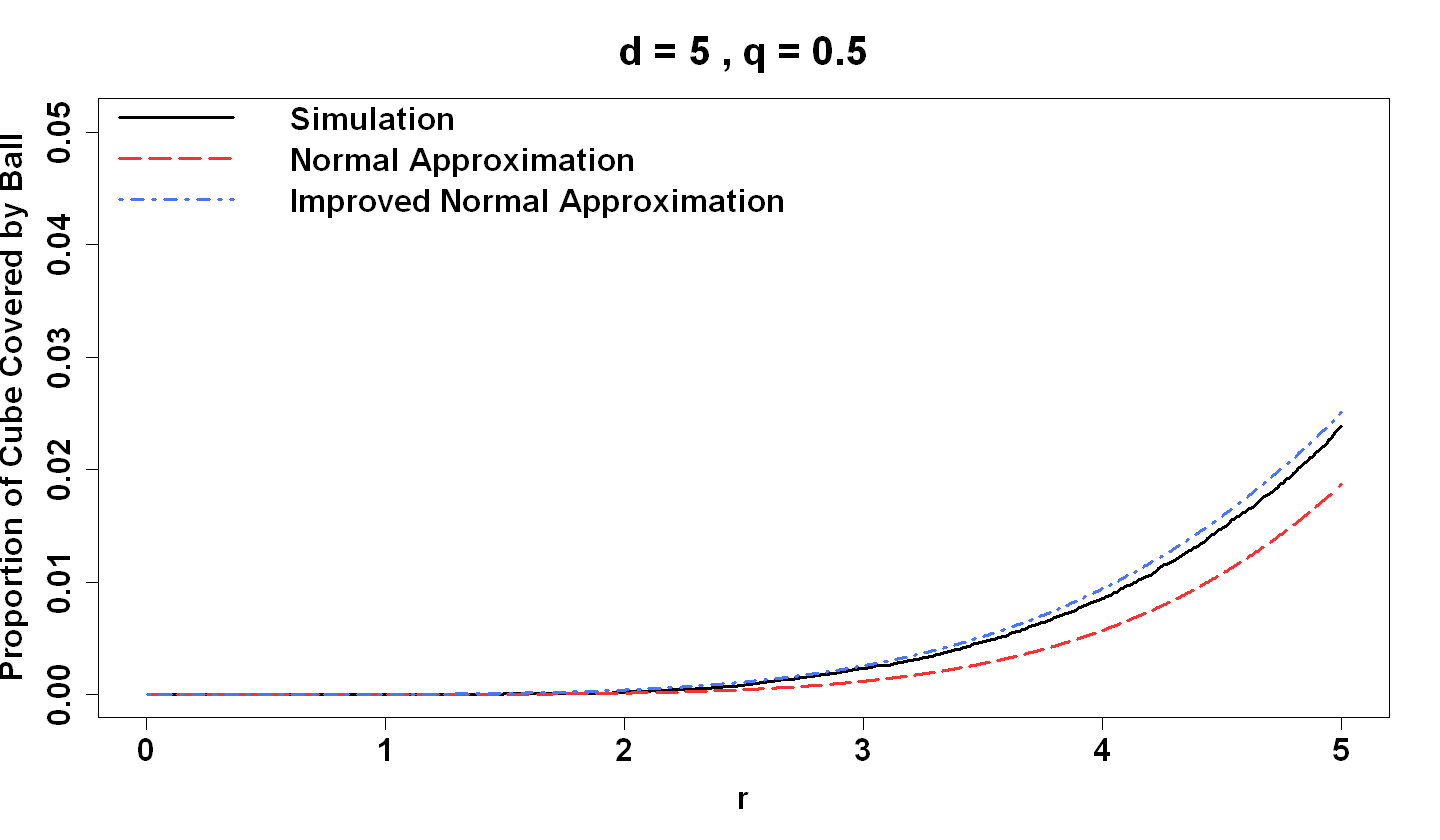}
  \captionsetup{justification=centering}
  \captionof{figure}{\textsl{Graphical representation of Monte Carlo simulation (black) for \(C_{d,q,r}\), \(q=0.5\), \(d=5\), with Normal (red) and improved Normal (blue) approximations.}}
  \label{fig:q0.5d5}
\end{minipage}
\begin{minipage}{.45\textwidth}
  \centering
  \includegraphics[width=0.9\linewidth]{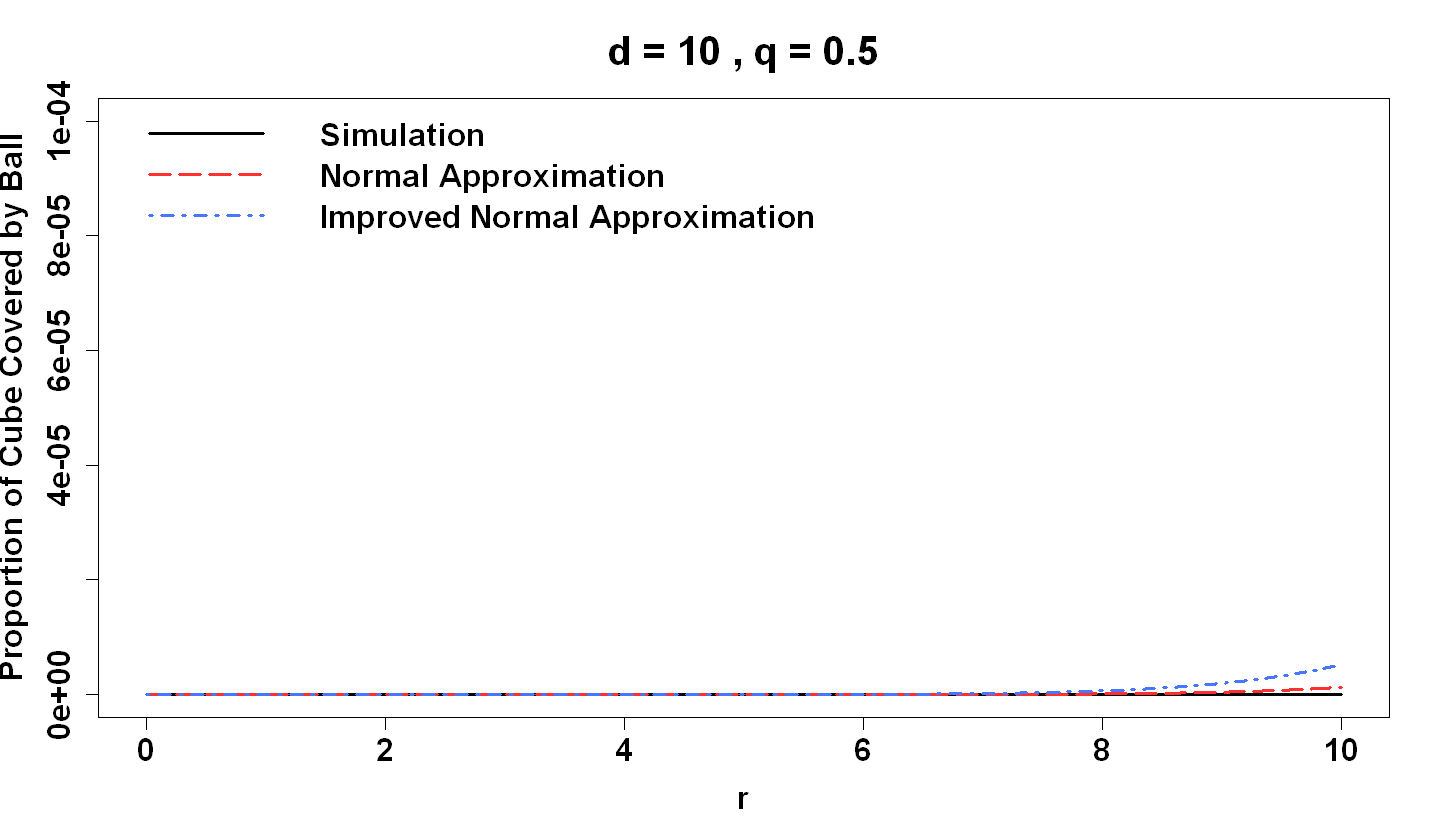}
  \captionsetup{justification=centering}
  \captionof{figure}{\textsl{Graphical representation of Monte Carlo simulation (black) for \(C_{d,q,r}\), \(q=0.5\), \(d=10\), with Normal (red) and improved Normal (blue) approximations.}}
  \label{fig:q0.5d10}
\end{minipage}
\end{figure}

In Figure \ref{fig:q0.5d5} we see that, unlike the other figures plotted over the whole range of possible \(r\), the Monte Carlo simulation (black) never reaches one. The improved Normal approximation (blue) overestimates very slightly but seems to follow the trend well. On the other hand the the first normal approximation (red) is not as accurate for \(r\) greater than around 3. We may expect though, that \(d=5\) would be too small for the approximations to hold in general.

Conversely, Figure \ref{fig:q0.5d10} shows that the simulation of the intersection volume remains at zero over the whole domain. On the one hand, this could possibly be due to not generating enough points in the cube. However, from Section 1.1.2 of \cite{AnatolyBook} we know that as the dimension increases, the mass of the cube is concentrated close to the boundary, whilst we may expect the shape of an \(L_{q}\)-ball for \(0 < q < 1\) to be significantly different to that for \(q \geq 1\). Thus, more simulations and more investigations into the theory for this case should be performed to gain a deeper understanding of this scenario.

\section{Conclusion}\label{sec:Conclusion}
In conclusion, this report has investigated the intersection volume between the cube \([-1,1]^{d}\) and an \(L_{q}\)-ball of radius \(r\) centred at the origin. This has been completed by describing two different approximation formulae for the quantity, exploring the depth of knowledge of ChatGPT in this area, and performing some simulations to compare the accuracy of the formulae for different values of \(q\) and \(d\).

Firstly, in Section \ref{sec:Theory}, we utilised the techniques of calculating moments, the Law of Large Numbers and the Central Limit Theorem to arrive at a Normal distribution based formula to approximate the fraction of the cube covered by the ball, given by equation (\ref{eq:Approx1BoxVersion}). We next used the first three cumulants and the first order Edgeworth Expansion to give a second, improved approximation formula based on the Normal distribution, which is stated in equations (\ref{eq:Approx2BoxVersion}) and (\ref{eq:Approx2BoxVersiontrValue}). This section could be extended further by considering including additional correction terms from the Edgeworth Expansion, or by broadening the scenario considered perhaps to an \(L_{q}\)-ball centred away from the origin or to multiple \(L_{q}\)-balls.

In Section \ref{sec:ChatGPT}, we asked ChatGPT a number of questions linked to the formulae presented in Section \ref{sec:Theory} and evaluated the responses obtained. We concluded that, whilst it attempts to provide some answers in this area, particularly in the more studied case of \(q = 2\), the information it returns is not completely reliable, with mistakes being made. To improve this section further, more questions should be asked to potentially test the reproducibility of the answers gained, or to see if alternative phrasings of the questions may lead to more accurate responses.

Finally, in Section \ref{sec:Simulations}, we performed many simulations to evaluate the accuracy of the approximation formulae for different values of \(q\) and \(d\). We found that, for \(q > 1\), increasing the dimension seems to improve the accuracy of both formulae, whilst the second one appears significantly more reliable than the first in almost all cases. Additionally, the dimension required to maintain the same level of accuracy seems to increase as \(q\) increases. Paying closer attention to the tails, it was observed that the second approximation formula can become negative, despite approximating a non-negative quantity. Furthermore, investigations suggested that the trend of values obtained in the case that \(0 < q < 1\) is greatly different to that for \(q \geq 1\). Thus, further experimentations into these cases should be performed to improve the report.

Overall, this report has explored the topic of approximating the volume of intersection of a cube and an \(L_{q}\)-ball in high dimensions by introducing some of the theory, exploring responses given by ChatGPT and producing graphs illustrating the formulae. It has also suggested some areas for further investigations that could improve and extend the report, enabling a more in depth understanding of the subject to be obtained.

\subsection*{Acknowledgement}
The author is very grateful to Professor Anatoly Zhigljavsky for attracting the author's attention to this problem and for providing useful advice.

\newpage
\section{References}
\small{
\printbibliography[heading=none]
}

\newpage
\appendix
\section{Appendix}\label{sec:Appendix}
In this section we present the code used to produce the figures seen in Section \ref{sec:Simulations}. We start by setting the parameters as shown in Section \ref{subsec:Parameters}. We change the values of \(d\) and \(q\) to produce the different figures. Also, when considering the graph of the tails, we increase \(N\) to increase the accuracy of the Monte Carlo simulation, but restrict the range of \(r\) considered in \emph{radii} to just the interval of interest to keep the run time reasonable. The code calculating the Monte Carlo simulation for the fraction of the cube covered by the ball for different values of \(r\) for the \(d\) and \(q\) fixed is shown in Section \ref{subsec:MonteCarlo}. Section~\ref{subsec:NormalApproxFormula} and Section~\ref{subsec:ImprovedNormalApproxFormula} respectively contain code computing the first and second approximation formulae. Finally, Section \ref{subsec:Plotting} shows an example of the code used to plot the graphs presented in Section \ref{sec:Simulations}.

\subsection{Code for Setting Parameters}\label{subsec:Parameters}
\begin{lstlisting}[breaklines=true, language=R, commentstyle=\color{teal}]
    # Change these values.
    d <- 25 # d = 5, 10, 25, 50.
    q <- 5  # q = 0.5, 1, 2, 5, 10.
    N <- 100000 # Increase to produce graphs for the tail.
    width <- 0.005
    radii <- seq(width,d,by=width) # Use narrower range for the tails.
\end{lstlisting}

\subsection{Code for Monte Carlo Simulation}\label{subsec:MonteCarlo}
\begin{lstlisting}[breaklines=true, language=R, commentstyle=\color{teal}]
    # Generate N uniformly distributed points in the cube [-1,1]^d.
    points_in_cube <- list()
    for (i in 1:N){
        points_in_cube[[i]] <- runif(d,min=-1,max=1)
    }

    # Find the qth power of the Lq norm of each point in points.
    qth_power_of_Lq_norm <- list()
    for (i in 1:N){
        qth_power_of_Lq_norm[[i]] <- sum(abs(points_in_cube[[i]])^q)
    }

    # Find the proportion of the points in the cube that are in the Lq ball of radius r for all r in radii.
    pb <- txtProgressBar(max=length(radii), style=3)

    proportion_in_ball <- list()
    for (r in 1:length(radii)){
        setTxtProgressBar(pb, r)
        number_in_ball <- 0
        for (i in 1:N){
            if (qth_power_of_Lq_norm[[i]] <= radii[r]^q) {
                number_in_ball <- number_in_ball + 1
            }
        }
        proportion_in_ball[[r]] <- number_in_ball / N
    }
\end{lstlisting}

\subsection{Code for Normal Approximation Formula}\label{subsec:NormalApproxFormula}
\begin{lstlisting}[breaklines=true, language=R, commentstyle=\color{teal}]
    # Find the arguments to pass to the CDF of the standard Normal distribution for each value of r in radii.
    approx_1_args <- (sqrt(2*q+1)*((q+1)*(radii^q)-d))/(q*sqrt(d))

    # Calculate the value of the first approximation formula at each argument in approx_1_args.
    approx_1 <- pnorm(approx_1_args)
\end{lstlisting}

\subsection{Code for Improved Normal Approximation Formula}\label{subsec:ImprovedNormalApproxFormula}
\begin{lstlisting}[breaklines=true, language=R, commentstyle=\color{teal}]
    # Find the arguments t_r to pass to the CDF of the standard Normal distribution, PDF of the standard Normal distribution and
    # (1-(t_r)^2) as defined in the second approximation formula.
    approx_2_args <- (sqrt(2*q+1)*((q+1)*(radii^q)-d))/(q*sqrt(d))

    # Find the value of the constant that multiplies the second term in the second approximation formula.
    const <- (2*(q-1)*sqrt(2*q+1))/(6*sqrt(d)*(3*q+1))

    # Calculate the value of the second approximation formula at each argument in approx_2_args.
    approx_2 <- pnorm(approx_2_args) + (1 - (approx_2_args^2))*dnorm(approx_2_args)*const
\end{lstlisting}

\subsection{Code for Plotting the Graphs}\label{subsec:Plotting}
\begin{lstlisting}[breaklines=true, language=R, commentstyle=\color{teal}]
    # Plot the Monte Carlo simulation of the fraction of cube covered by the ball, the first approximation formula and the second
    # approximation formula over the full range of r in radii. Change xlim, ylim and the position of the legend for the different values of d and q used.
    options(repr.plot.width=12, repr.plot.height=7)

    plot(radii,proportion_in_ball,type='l',xlim=c(0,5),ylim=c(0,1),lty=1,lwd=3,col='black',
         xlab='r',ylab='Proportion of Cube Covered by Ball',xaxt='n',yaxt='n',cex.lab=1.6,font.lab=2,
         main=paste('d =',d, ', q =', q),cex.main=2,font.main=2)
    axis(1,cex.axis=1.6,font=2)
    axis(2,cex.axis=1.6,font=2)
    lines(radii,approx_1,type='l',col='firebrick1',lty=5,lwd=3)
    lines(radii,approx_2,type='l',col='royalblue1',lty=4,lwd=3)
    legend(2,0.25, legend=c('Simulation','Normal Approximation','Improved Normal Approximation'),
           col=c('black','firebrick1','royalblue1'), lwd=c(3,3,3), lty=c(1,5,4), horiz=FALSE, bty='n',text.font=2,cex=1.6)
\end{lstlisting}

\end{document}